\documentstyle[12pt, macro, amssymb, hyperref]{amsart}

\newtheorem{thm}{Theorem}[section]
\newtheorem{cor}[thm]{Corollary}
\newtheorem{prop}[thm]{Proposition}
\newtheorem{clm}[thm]{Claim}
\newtheorem{lem}[thm]{Lemma}

\theoremstyle{definition}
\newtheorem*{defn}{Definition}

\newtheorem*{rem}{Remark}
\newtheorem*{rems}{Remarks}

\theoremstyle{remark}

\numberwithin{equation}{section}

\def\subs{\subseteq}
\def\normal{\triangleleft}
\def\iso{\cong}
\def\iff{{\quad\Leftrightarrow\quad}}
\def\symdif{\triangle}
\def\bs{{\backslash}}

\newcommand{\ol}[1]{\overline{#1}}
\newcommand{\ben}{\begin{enumerate}}
\newcommand{\een}{\end{enumerate}}
\newcommand{\beit}{\begin{itemize}}
\newcommand{\enit}{\end{itemize}}
\newcommand{\bqr}{\begin{eqnarray}}
\newcommand{\eqr}{\end{eqnarray}}
\newcommand{\beq}{\begin{equation}}
\newcommand{\eeq}{\end{equation}}
\newcommand{\bqrn}{\begin{eqnarray*}}
\newcommand{\eqrn}{\end{eqnarray*}}

\def\pprime{{\prime\prime}}
\def\qed{\hfill$\framebox{}$\break\ \break}

\def\subs{\subseteq}
\def\GL{{\rm GL}}
\def\SL{{\rm SL}}
\def\PSL{{\rm PSL}}

\def\SO{{\rm SO}}
\def\SU{{\rm SU}}

\def\Rel{{\sf R}}

\def\Pf{\noindent{\bf Proof.\ }}

\def\Id{{\rm Id}}
\def\Ad{{\rm Ad}\,}

\def\rk{{\rm rk}_{\bR}}

\def\sB{\mathcal{B}}

\def\rmand{{\qquad{\rm and}\qquad}}

\def\rmwith{{\qquad{\rm with}\qquad}}
\def\rmwhere{{\qquad{\rm where}\qquad}}

\def\Gaut{{\rm Aut}\,}

\def\Gout{{\rm Out}\,}
\def\Raut{{\rm Aut}\,}
\def\Rinn{{\rm Inn}\,}
\def\Rout{{\rm Out}\,}
\def\Maut{{\rm Aut}}
\def\Aaut{{\rm Aut}}
\def\OA{{\sf A}}
\def\Aff{{\rm Aff}}

\title[Outer Automorphism Groups]{Outer~Automorphism groups of some
Ergodic~Equivalence~Relations}

\author{Alex Furman}\thanks{Supported in part by NSF
grants DMS-0049069 and DMS-0094245.}

\address{Dept of Mathematics (m/c 249),
University of Illinois at Chicago, 851 S. Morgan Street, Chicago,
IL 60607.}
\email{furman@@math.uic.edu}


\begin{document}

\begin{abstract}
Let $\Rel$ a be countable ergodic equivalence relation of type
${\rm II}_1$ on a standard probability space $(X,\mu)$. The group
$\Rout\Rel$ of \emph{outer automorphisms} of $\Rel$ consists of
all invertible Borel measure preserving maps of the space which
map $\Rel$-classes to $\Rel$-classes modulo those which preserve
almost every $\Rel$-class. We analyze the group $\Rout\Rel$ for
relations $\Rel$ generated by actions of higher rank lattices,
providing general conditions on finiteness and triviality of
$\Rout\Rel$ and explicitly computing $\Rout\Rel$ for the standard actions.
The method is based on Zimmer's superrigidity for measurable cocycles,
Ratner's theorem and Gromov's Measure Equivalence construction.
\end{abstract}

\maketitle

\section{Introduction and Statement of the Main Results}

Let $(X,\sB)$ be a standard Borel space with a non-atomic
probability Lebesgue measure and let $\Rel$ be a countable
measurable relation of type ${\rm II}_1$ on $(X,\sB,\mu)$, i.e.
measurable, countable, ergodic and measure preserving equivalence
relation $\Rel\subset X\times X$. For the abstract definition of
this notion the reader is referred to the fundamental work of
Feldman and Moore \cite{FM}, which in particular demonstrates that
any such equivalence relation can be presented as the orbit
relation
\[
   \Rel_{X,\GA} = \{ (x,y)\in X\times X \mid \GA\cdot x = \GA\cdot y \}
\]
of an \emph{ergodic, measure preserving} action of some countable
group $\GA$ on the space $(X,\sB,\mu)$. In most of the examples in
this paper equivalence relations are defined by ergodic
measure-preserving actions of concrete countable groups $\GA$,
namely irreducible lattices in semi-simple
connected higher rank real Lie groups.

In the purely measure-theoretical context of this paper all
objects are considered modulo sets of zero $\mu$-measure, denoted (mod 0).
Following this convention the measure space automorphism
group $\Maut(X,\mu)$ is the group of all invertible Borel maps
$T:X\to X$ with $T_*\mu=\mu$, where two such maps which
agree on a set of full $\mu$-measure are identified. In a similar
fashion two equivalence relations $\Rel$, $\Rel^\prime$ on
$(X,\mu)$ are identified if there exists a subset $X^\prime\subs X$
with $\mu(X^\prime)=1$ on which the restrictions of $\Rel$ and
$\Rel^\prime$ coincide.

Given an equivalence relation $\Rel$ on $(X,\mu)$ consider the group of
\emph{relation automorphisms}
\[
   \Raut\Rel=\{T\in\Maut(X,\mu) \mid T\times T(\Rel)=\Rel \}
\]
and the subgroup $\Rinn\Rel$ of \emph{inner automorphisms},
also known as the \emph{full group} of $\Rel$,
consisting of such $T\in\Maut(X,\mu)$ that $(x,Tx)\in\Rel$ for
$\mu$-a.e. $x\in X$. The full group $\Rinn\Rel$ is normal in
$\Raut\Rel$ and the \emph{outer automorphism group}
$\Rout\Rel$ is defined as the quotient
\[
   1\overto{}\Rinn\Rel\overto{}\Raut\Rel
   \overto{\ep}\Rout\Rel\overto{}1
\]
Elements of $\Rout\Rel$ represent measurable ways to permute
$\Rel$-classes on $(X,\mu)$. The full group $\Rinn\Rel$ is always
very large (see Lemma~\ref{L:full}). For the unique
\emph{amenable} equivalence relation $\Rel_{am}$ of type ${\rm
II}_1$ the outer automorphism group $\Rout\Rel_{am}$ is also
enormous. The purpose of this paper is to analyze
$\Rout\Rel_{X,\GA}$ for orbit relations $\Rel_{X,\GA}$ generated
by m.p. ergodic actions of higher rank lattices, in particular
presenting many natural examples of relations $\Rel$ with trivial
$\Rout\Rel$. Such examples were first constructed by S.~Gefter in
\cite{Ge1}, \cite{Ge2} (Thm~\ref{T:ex-compact} below).

Prior to stating the results let us define two special subgroups
in $\Rout\Rel$, in the case where $\Rel$ is the orbit relation
$\Rel_{X,\GA}$ generated by some measure-preserving action
$(X,\mu,\GA)$ of some countable group $\GA$. In such a situation
consider the group $\Aaut(X,\GA)$ of \emph{action automorphisms}
of the system $(X,\mu,\GA)$
\[
   \Aaut(X,\GA):=\{T\in\Maut(X,\mu) \mid
      T(\ga\cdot x)=\ga\cdot T(x),\ \forall\ga\in\GA\}.
\]
This is the centralizer of $\GA$ in $\Maut(X,\mu)$. For a group
automorphism $\ta\in\Gaut\GA$ define
\[
   \Aaut^\ta(X,\GA):=
      \left\{T\in\Maut(X,\mu) \mid T(\ga\cdot x)=\ga^\ta\cdot T(x)
      \right\}
\]
and let $\Aaut^*(X,\GA)$ be the union of $\Aaut^\ta(X,\GA)$ over
$\ta\in\Gaut\GA$.  (If the $\GA$-action is faithful
$\Aaut^*(X,\GA)$ is the normalizer of $\GA$ in $\Maut(X,\mu)$). We
shall denote by $\OA(X,\GA)$ and $\OA^*(X,\GA)$ the images of the
groups $\Aaut(X,\GA)$ and $\Aaut^*(X,\GA)$ under the factor map
$\Raut\Rel_{X,\GA}\overto{\ep}\Rout\Rel_{X,\GA}$. Observe that the
$\ep$-image in $\Rout\Rel_{X,\GA}$ of the coset $\Aaut^\ta(X,\GA)$
depends only on the outer class $[\ta]\in\Gout\GA$ and therefore
can be denoted by $\OA^{[\ta]}(X,\GA)$. The group $\OA(X,\GA)$ is
normal in $\OA^*(X,\GA)$ and the factor group
$\OA^*(X,\GA)/\OA(X,\GA)$ is (a factor of) a subgroup
of $\Gout\GA$. 
In general, the subgroups
$\OA(X,\GA)\subs\OA^*(X,\GA)$ of $\Rout\Rel_{X,\GA}$ depend on the
specific presentation of the relation $\Rel$ as the orbit relation
$\Rel_{X,\GA}$ of an action $(X,\mu,\GA)$.

In this paper we are mostly interested in ergodic m.p.
actions of higher rank lattices and will be using the following terminology
and notations:
\begin{itemize}
\item
  For locally compact, secondly countable group $G$ a left-invariant
  Haar measure will be denoted by $m_G$. If $\GA\subset G$ is
  a discrete group so that $G/\GA$ carries a finite $G$-invariant measure
  we say that $\GA$ forms a \emph{lattice} in $G$ and will denote by
  $m_{G/\GA}$ the unique $G$-invariant probability measure on $G/\GA$.
\item
  The term \emph{semi-simple Lie group} will mean
  semi-simple, connected, center-free, real Lie group $G=\prod G_i$ without
  non-trivial compact factors, unless stated otherwise.
  A lattice $\GA$ in a semi-simple Lie group $G=\prod G_i$
  is \emph{irreducible} if $\GA$ does not contain a finite index subgroup
  $\GA^\prime$ which splits as a direct product of lattices in subfactors.
  By \emph{higher rank lattice} hereafter we shall mean an irreducible
  lattice in a semi-simple Lie group $G$ with $\rk(G)\ge 2$.
\item
  A measure-preserving action $(X,\mu,\GA)$ of a lattice $\GA$ in a
  semi-simple Lie group $G=\prod G_i$ is \emph{irreducible} if the action
  of every simple factor $G_i$ in the induced
  $G$-action on $(G\times_\GA X,m_{G/\GA}\times\mu)$ is ergodic.
  Clearly, if $G$ is \emph{simple} then any lattice $\GA\subset G$ is irreducible
  and any ergodic $\GA$-action is irreducible.
\item
  For an arbitrary group $\GA$ a m.p. action $(X_0,\mu_0,\GA)$ is a
  ($\GA$-equivariant) \emph{quotient} of another m.p. action $(X,\mu,\GA)$
  if there exists a measurable map $\pi:X\to X_0$ with $\pi_*\mu=\mu_0$
  and $\pi(\ga\cdot x)=\ga\cdot\pi(x)$ for $\mu$-a.e. $x\in X$ and all
  $\ga\in\GA$.
\item
  A measure-preserving action $(X,\mu,\GA)$ of an arbitrary group $\GA$
  is called \emph{aperiodic} if every finite atomic quotient of
  $(X,\mu,\GA)$ is trivial; equivalently if every finite index subgroup
  $\GA^\prime\subset \GA$ acts ergodically on $(X,\mu)$.
\end{itemize}
\begin{rems}
\begin{itemize}
\item[(a)]
  Every \emph{mixing} ergodic action
  $(X,\mu,\GA)$ of an irreducible lattice $\GA$ in a semi-simple Lie
  group $G$ is irreducible and aperiodic.
\item[(b)]
  By the result of Stuck and Zimmer \cite{SZ} any ergodic non-atomic m.p.
  action of an irreducible lattice $\GA$ in a semi-simple Lie group $G$
  with property (T) is \emph{free} (mod 0). Recall that a higher rank
  semi-simple $G$ has property (T) iff it does not have simple factors
  locally isomorphic to $\SO(n,1)$ or $\SU(n,1)$.
\item[(c)]
  For any free, ergodic action $(X,\mu,\GA)$ of an irreducible
  lattice $\GA$ in a semi-simple Lie group $G$ the map
  \[
      \Aaut(X,\GA)\overto{\ep}\OA(X,\GA)
  \]
  is an isomorphism and the homomorphism
  \[
      \OA^*(X,\GA)/\OA(X,\GA)\overto{}\Gout\GA
  \]
  is an embedding (Lemma~\ref{L:ICC} below).
\item[(d)]
  It follows from the Strong Rigidity (Mostow, Prasad, Margulis)
  that for an irreducible lattice
  $\GA\subset G\not\simeq\SL_2(\bR)$ the automorphism group $\Gaut\GA$ is
  isomorphic to the normalizer
  $N_{\Gaut G}(\GA)$ of $\GA$ in $\Gaut G\supseteq \Ad G\iso G$. Since
  $\GA_*:=N_{\Gaut G }(\GA)\supseteq\GA$ is a closed subgroup
  properly contained
  in $\Gaut G$, it forms a lattice in $\Gaut G$, and
  $\Gout\GA =\GA_*/\GA$ is always finite.
\end{itemize}
\end{rems}

Thus for an irreducible aperiodic free m.p. action of a higher
rank lattice $\GA$ the analysis of $\Rout\Rel_{X,\GA}$ reduces to
the analysis of the quotient $\Rout\Rel_{X,\GA}/\OA^*(X,\GA)$ and
the subgroup $\OA^*(X,\GA)$ which, up to at most finite index, is
isomorphic to $\Aaut(X,\GA)$.
\begin{thm} \label{T:main}
Let $G$ be a semi-simple, connected, center-free, real Lie group
without non-trivial compact factors and with $\rk(G)\ge 2$. Let
$\GA\subset G$ be an irreducible lattice and $(X,\mu,\GA)$ be a
measure preserving, ergodic, irreducible, aperiodic, essentially
free $\GA$-action. Assume that $(X,\mu,\GA)$ does not admit
measurable $\GA$-equivariant quotients of the form
$(G/\GA^\prime,m_{G/\GA^\prime},\GA)$ where $\GA^\prime\subset G$
is a lattice isomorphic to $\GA$ and $\GA$ acts by
$\ga:g\GA^\prime\mapsto\ga g\GA^\prime$. Then
\[
    \Rout\Rel_{X,\GA}=\OA^*(X,\GA)
\]
while $\OA^*(X,\GA)\iso\Aaut(X,\GA)/\GA$.
\end{thm}
More generally
\begin{thm} \label{T:mainN}
Let $\GA\subset G$ be a higher rank lattice as in
Theorem~\ref{T:main} and $(X,\mu,\GA)$ be any measure preserving,
ergodic, irreducible, aperiodic, essentially free $\GA$-action.
If  $OA^*(X,\GA)$ has finite index $n>1$ in $\Rout\Rel_{X,\GA}$ 
then $(X,\mu,\GA)$ has an equivariant measurable quotient
\[
     \pi\,:\, (X,\mu)\overto{} (G^{n-1}/\GA^{n-1},m_{G^{n-1}/\GA^{n-1}})
                               =\prod_{i=1}^{n-1} (G/\GA,m_{G/\GA})
\]
where the $\GA$-action on
$(G^{n-1}/\GA^{n-1},m_{G^{n-1}/\GA^{n-1}})$ is given by
\[
     \ga\,:\,(x_i)_{i=1}^{n-1}\,\mapsto\,(\ga^{\ta_i}\cdot x_i)_{i=1}^{n-1}
\]
for some fixed automorphisms $\ta_i\in\Gaut G$, $1\le i<n$.

If $\OA^*(X,\GA)]$ has infinite index in $\Rout\Rel_{X,\GA}$ 
then $(X,\mu,\GA)$ has an 
infinite product equivariant quotient space
\[
     \pi\,:\, (X,\mu)\overto{} \prod_{i=1}^{\infty} (G/\GA,m_{G/\GA})
\]
with a diagonal $\GA$-action on 
$\ga\,:\,(x_i)_{i=1}^{\infty}\mapsto (\ga^{\ta_i}x_i)_{i=1}^\infty$
for some fixed sequence $\ta_i\in\Gaut G$, $i=1,2,\dots$.
\end{thm}

Of course, Theorem~\ref{T:main} is just a particular case of \ref{T:mainN}
(contrapositive formulation for $n=1$) since a $\ta$-twisted
$\GA$-action $\ga:g\GA\mapsto\ga^\ta g\GA$ on $(G/\GA,m_{G/\GA})$
is measurably isomorphic to the untwisted $\GA$-action
$\ga:g\GA^\prime\mapsto \ga g\GA^\prime$ where
$\GA^\prime=\ta^{-1}(\GA)$.

For an $d\times d$ matrix $A$ let $\ch(A):=\sum_i
\log^+|\la_i(A)|$, where $\log^+x=\max\{0,\log{x}\}$ and
$\la_i(A)$ denote the eigenvalues of $A$. Given a semi-simple
group $G$ and $d\in\bN$ consider all linear representations
$\rh:G\to\GL_d(\bC)$ (there are finitely many equivalence classes
for any $d$) and let
\[
   W_G(d):=\max_{\dim\rh=d}
   \inf_{g\in G}\ \frac{\ch(\rh(g))}{\ch(\Ad(g))}
\]
\begin{cor}\label{C:entropy}
Let $\GA\subset G$ and $(X,\mu,\GA)$ be as in
Theorem~\ref{T:mainN}. Denote by $h(X,\ga)$ the Kolmogorov-Sinai
entropy of the single measure-preserving transformation $\ga$ of
$(X,\mu)$. Then
\begin{equation}\label{e:entropy}
   [\Rout\Rel_{X,\GA}:\OA^*(X,\GA)]\le
   1+\inf_{\ga\in\GA}\ \frac{h(X,\ga)}{\ch(\Ad(\ga))}
\end{equation}
If $X$ is a compact manifold with a $C^1$-action of a higher rank lattice
$\GA\subset G$ which preserves a probability measure $\mu$ on $X$
so that $(X,\mu,\GA)$ is a free (mod 0) action which is
ergodic, irreducible and aperiodic, then
\begin{equation}\label{e:estimate}
    [\Rout\Rel_{X,\GA}:\OA^*(X,\GA)]\le 1+W_G(\dim(X))
\end{equation}
The function $W_G$ satisfies $W_G(d)\le d^2/8$.
\end{cor}
\begin{rem}
Theorem~\ref{T:G-mod-GA} below shows that the inequality
(\ref{e:entropy}) is sharp. However the estimate
(\ref{e:estimate}) is probably not optimal, with a more plausible
one being $1+\dim(X)/\dim\,{\rm Lie}(G)$.
\end{rem}
\begin{rem}
As we shall see below, groups $\Rout\Rel_{X,\GA}$ and $\OA(X,\GA)$ 
can be very large when considered as abstract groups, but in all cases 
below the quotient $\Rout\Rel_{X,\GA}/\OA^*(X,\GA)$ is either finite or countable.
This might be a general property of actions of higher rank lattices. 
In fact, this property is known for essentially free ergodic actions $(X,\GA)$ of groups $\GA$ 
with Kazhdan's property (T).
For such groups (and in a slightly more general situation) Gefter and Golodets introduced a 
natural topology on $\Rout\Rel_{X,\GA}$ with respect to which 
$\Rout\Rel_{X,\GA}$ is a Polish (i.e. complete separable) group
and $\OA(X,\GA)$ is an open subgroup
(see \cite{GeG} Thm 2.3 and references throughout section 2). 
\end{rem}

In specific cases, in particular in the standard examples of
algebraic lattice actions, it is possible to compute the groups 
$\Rout\Rel_{X,\GA}$ explicitly as we shall do in 
Theorems~\ref{T:torus} - \ref{T:G-mod-GA} below. 
In Theorems~\ref{T:torus} - \ref{C:ex-H-LA} the systems 
$(X,\GA)$ do not have $G/\GA^\prime$ as measurable 
quotients and therefore by Theorem~\ref{T:main} we have 
$\Rout\Rel_{X,\GA}=\OA^*(X,\GA)=\Aaut^*(X.\GA)/\GA$.
The latter groups are of algebraic nature, but their explicit descriptions
are cumbersome. Thus for readers convenience we have also 
presented the groups $\Aaut(X,\GA)$, which have a more transparent appearance
and have at most finite ($\le|\Gout\GA|$)
index in $\OA^*(X,\GA)$.
\begin{thm} \label{T:torus}
Let $G$ be a simple, connected real Lie group with finite center
and $\rk(G)\ge 2$ and $\rh:G\hookrightarrow\SL_N(\bR)$ be an
embedding such that $\rh(G)$ does not have non-trivial fixed
vectors and assume that $G$ has a lattice $\GA\subset G$ such that
$\rh(\GA)\subs\SL_N(\bZ)$. Then the natural $\GA$-action on the
torus $\bT^N=\bR^N/\bZ^N$ is ergodic, aperiodic and the orbit
relation $\Rel_{\bT^N,\GA}$ satisfies
\[
  \begin{array}{rllll}
      \Rout\Rel_{\bT^N,\,\GA} &=& \OA^*(\bT^N,\GA) &\iso&
       N_{\GL_N(\bZ)}(\rh(\GA))/\rh(\GA) \\
      \OA(\bT^N,\GA)&\iso&\Aaut(\bT^N,\GA) &\iso&
      C_{\GL_N(\bR)}(\rh(G))\cap\GL_N(\bZ)
  \end{array}
\]
In particular, for $n>2$ the $\SL_n(\bZ)$-action on $\bT^n$ gives
an ergodic relation $\Rel_{\bT^n,\SL_n(\bZ)}$ which has no outer
automorphisms if $n$ is even, and a single outer automorphism
given by $x\mapsto -x$ if $n$ is odd.
\end{thm}
Note that in the above Theorem we allowed finite non-trivial centers to
accommodate the standard example of $\GA=\SL_n(\bZ)$ acting on the
torus $\bT^n$ for even $n>2$.

To state the following results we recall the notion of
\emph{affine transformations} of a homogeneous space
(these are needed only for the precise description of $\Rout\Rel_{X,\GA}$),
however the spirit of the results is captured by the finite index subgroup
$\OA(X,\GA)$ which does not require this notion.
\begin{defn}\label{D:affine}
Let $\LA$ be a subgroup of a group $H$, and let $N:=\bigcap_{h\in
H} h^{-1}\LA h$ denote the maximal subgroup of $\LA$ which is
normal in $H$. Given an automorphism $\si$ of $H/N$ with
$\si(\LA/N)=\LA/N$ and $t\in H/N$ denote by $a_{\si,t}$ the map
\[
   a_{\si,t}:h\LA\mapsto t\si(h)\LA
\]
of $H/\LA$. Such maps will be called \emph{affine}, and we shall
denote by $\Aff(H/\LA)$ the group of all affine maps of $H/\LA$.
\end{defn}
Replacing $H$ by $H/N$ and $\LA$ by $\LA/N$ one does not change
the homogeneous space: $H/\LA\iso(H/N)/(\LA/N)$. Thus hereafter we
shall assume that $N$ is trivial. Under this assumption the map
$(\si,t)\mapsto a_{\si,t}$ defines an epimorphism $
  N_{\Gaut H}(\LA)\ltimes H \overto{} \Aff(H/\LA)
$
(which contains $\{(\Ad\la,\la)\mid \la\in\LA\}$ in its kernel)
and which maps $H\iso\{\Id\}\times H$ isomorphically onto its
image in $\Aff(H/\LA)$. This copy of $H$ in $\Aff(H/\LA)$ has
index bounded by $|\Gout H|$. 

We shall be interested in situations where some group (a higher rank 
lattice) $\GA$ is embedded in $H$, $\rh:\GA\to H$, and acts on the
homogeneous space $H/\LA$ by left translations.
Then the normalizer $N_{\Aff(H/\LA)}(\rh(\GA))$ in $\Aff(H/\LA)$
of this action consists of those affine maps
$a_{\si,t}$ for which $\si\in\Gaut H$ and $t\in H$ satisfy
\[
   \si(\LA)=\LA\rmand \si(\rh(\GA))=t^{-1}\rh(\GA)t.
\]
In any case this group contains $N_H(\rh(\GA))$ as a subgroup
of an index bounded by $|\Gout H|$.
\begin{thm}[cf. Gefter \cite{Ge2}]\label{T:ex-compact}
Let $\GA$ be a higher rank lattice which admits a dense embedding
$\rh:\GA\to K$ into a compact connected Lie group $K$. Then for
every closed subgroup $\{e\}\subs L\subset K$ the $\GA$-action on
$(K/L,m_{K/L})$ is ergodic, irreducible and aperiodic and the
orbit relation $\Rel_{K/L,\GA}$ satisfies
\[
  \begin{array}{rllll}
    \Rout \Rel_{(K/L,\,\GA)}&=&\OA^*(K/L,\GA) &\iso& N_{\Aff(K/L)}(\rh(\GA))/\rh(\GA)\\
    \OA(K/L,\GA)&\iso&\Aaut(K/L,\GA)&\iso& N_K(L)/L
  \end{array}
\]
In particular, if the compact group $K$ has no outer automorphisms 
which normalize $L$ or if $\GA$ has no outer automorphisms, then
\[
	\Rout\Rel_{K/L,\GA}\iso N_K(L)/L.
\]
\end{thm}
\begin{rem}
  Theorem~\ref{T:ex-compact} was proved by S.~Gefter in \cite{Ge2}, where he  
  constructed the first and only example of type ${\rm II}_1$ equivalence relations without outer
  automorphisms. Indeed, by a well known arithmetic construction (cf. \cite{Zi-book}
  5.2.12) certain lattices $\GA\subset G:=\SO(p,q)$ admit dense
  embeddings into the compact group $K:=\SO(n)$ where $n=p+q$. Take
  $p>q\ge 2$ to ensure $\rk(G)\ge 2$ and let $L\iso\SO(n-1)$ be the
  stabilizer of a point in the  in $K=\SO(n)$ action on the sphere $S^{n-1}$. 
  Then $N_K(L)=L$ and
  since $\SO(n)$ has no outer automorphisms, $\Aff(K/L)\iso
  N_K(L)$, which shows that $\Rout\Rel_{(K/L,\GA)}$ is trivial.
\end{rem}
In Theorem~\ref{T:ex-compact} the compact group $K$ is taken to be
a \emph{connected} Lie group to guarantee aperiodicity of the
action. Higher rank lattices can also be
densely embedded in other compact groups, namely profinite ones. 
Such embeddings give rise to ergodic actions which strongly violate
aperiodicity condition - they are inverse limits of finite quotients.
A typical example is the standard embedding
\[
  \GA:=\PSL_n(\bZ) \overto{\rh} K:=\PSL_n(\bZ_p).
\]
It was observed by Gefter (\cite{Ge2} Rem 2.8) that in this case
$\Rout\Rel_{K,\GA}$ contains a group isomorphic to $\PSL_n(\bQ_p)$,
in such a way that 
\[
	\OA(K,\GA)\iso K=\PSL_n(\bZ_p)\subset \PSL_n(\bQ_p)\subs \Rout\Rel_{K,\GA}
\]  
so that $\OA(K,\GA)$ has infinite index in $\Rout\Rel_{K,\GA}$.
We claim that the last inclusion is essentially an equality.
More generally
\begin{thm}\label{T:SLnZp}
Consider the natural dense embedding of $\GA=\PSL_n(\bZ)$, $n\ge 3$, 
in the profinite group 
$K=\prod_{i=1}^r \PSL_n(\bZ_{p_i})$ where $\{p_1,\dots,p_r\}$ is
a finite set of distinct primes. Then  
$\Rout\Rel_{K,\GA}$ is a $\bZ/2$ extension of
\[
	H=\prod_{i=1}^r\PSL_n(\bQ_{p_i})
\]
with the transpose map 
$(k_1,\dots,k_r)\mapsto (k_1^t,\dots,k_r^t)$ of $K$ giving rise to the $\bZ/2$ extension.
\end{thm}
Another family of standard examples is described by the following
\begin{thm}\label{T:ex-homogeneous}
Let $\GA\subset G$ be a higher rank lattice as in
Theorem~\ref{T:main}, $H$ be a connected Lie group, $\LA\subset H$
be a closed subgroup so that $H/\LA$ carries an $H$-invariant
probability measure $m_{H/\LA}$, and assume that $H$ does not
admit surjective homomorphisms $\si:H\to G$ with
$\si(\LA)\subs\GA$. Suppose that there exists a homomorphism
$\rh:G\to H$ such that each of the simple factors $G_i$ of $G$
acts ergodically on $(H/\LA,m_{H/\LA})$. Then for the $\GA$-action
on $(H/\LA,m_{H/\LA})$ one has
\[
  \begin{array}{rllll}
     \Rout\Rel_{(H/\LA,\,\GA)}&=& \OA^*(H/\LA,\,\GA)
         &\iso & N_{\Aff(H/\LA)}(\rh(\GA))/\rh(\GA)\\
     \OA(H/\LA,\GA)&\iso&\Aaut(H/\LA,\GA) &\iso& C_{\Aff(H/\LA)}(\rh(\GA))
  \end{array}
  \]
\end{thm}
\begin{cor}\label{C:ex-H-LA}
Let $\GA\subset G$ be a higher rank lattice as in
Theorem~\ref{T:main}, $H$ be a connected semi-simple Lie group
with trivial center, $\rh:G\hookrightarrow H$ be an embedding and
let $\LA\subset H$ be an irreducible lattice. Assume that either
 $\rh$ is a proper embedding, i.e. $G\not\simeq H$, or that
 $\rh:G\to H$ is an isomorphism but $\LA$ is not abstractly isomorphic to a
 subgroup of finite index in $\GA$.
Then the $\GA$-action on $(H/\LA,m_{H/\LA})$ by left translations
is ergodic, irreducible and aperiodic and the orbit relation
$\Rel_{H/\LA,\GA}$ has
\[
     \Rout\Rel_{(H/\LA,\GA)}=\OA^*(H/\LA,\GA)\iso
     N_{\Aff(H/\LA)}(\rh(\GA))/\rh(\GA)
\]
This group contains the centralizer $C_{H}(\rh(G))$ as a normal
subgroup of finite index dividing $|\Gout\LA|\cdot|\Gout\GA|$.
\end{cor}
\begin{rem}
Corollary~\ref{C:ex-H-LA} also allows to construct ergodic
equivalence relations without outer automorphisms. Indeed if a
simple Lie group $G\not\simeq\SL_2(\bR)$ has no outer
automorphisms, then maximal lattices $\GA$ in $G$ have trivial
$\Gout\GA$ as well. Choosing two non-commensurable maximal
lattices $\GA,\LA$ in such a $G$ one obtains an equivalence
relation $\Rel_{G/\LA,\GA}$ without outer automorphisms.
Similarly, one can find proper embeddings $G\subset H$ where $G$
and $H$ are simple higher rank Lie groups with $\Gout G$, $\Gout
H$, $C_H(G)$ all being trivial. Then for any choice of maximal
lattices $\GA\subset G$, $\LA\subset H$, the $\GA$-action on
$H/\LA$ gives $\Rel_{H/\LA,\GA}$ without outer automorphisms.
\end{rem}
All the examples discussed so far had the property that the
original system $(X,\mu,\GA)$ did not admit measurable
$\GA$-equivariant quotients of the form
$(G/\GA^\prime,m_{G/\GA^\prime},\GA)$; and therefore
Theorem~\ref{T:main} allowed to conclude that
$\Rout\Rel_{X,\GA}=\OA^*(X,\GA)\iso\Aaut^*(X,\GA)/\GA$. The
following result analyzes what happens if this assumption is not
satisfied.
\begin{thm}\label{T:G-mod-GA}
Let $\GA\subset G$ be a higher rank lattice as in
Theorem~\ref{T:main} and let $\GA$-act on $(G/\GA,m_{G/\GA})$ by
left translations. Then for the corresponding orbit relation
$\Rel_{G/\GA,\GA}$ 
  \bqrn
     [\Rout\Rel_{G/\GA,\GA}&:&\OA^*(G/\GA,\GA)]=2\\
     \OA(G/\GA,\GA)&\iso&\Aaut(G/\GA,\GA) \iso  \{1\}\\
     \OA^*(G/\GA,\GA) &\iso& \Gout\GA \\
     \Rout\Rel_{(G/\GA,\,\GA)} &\iso& (\bZ/2\bZ)\times \Gout\GA
  \eqrn
More generally, for any $n\in\bN$ the diagonal
left $\GA$-action on the product space
$(G^n/\GA^n,\,m_{G^n/\GA^n})$ satisfies
\[
  \begin{array}{rllll}
     [\Rout\Rel_{(G^{n}/\GA^{n},\,\GA)}&:&\OA^*(G^{n}/\GA^{n},\GA)] \ =\ n+1\\
     \OA(G^n/\GA^n,\GA)&\iso&\Aaut(G^{n}/\GA^{n},\GA) \ \iso\ S_{n}\\
     \OA^*(G^{n}/\GA^{n},\GA) &\iso& S_{n}\times \Gout\GA \\
     \Rout \Rel_{(G^{n}/\GA^{n},\,\GA)} &\iso& S_{n+1}\times
     \Gout\GA
  \end{array}
\]
where $S_n$ denotes the permutation group on
$\{1,\dots,n\}$.

For the diagonal $\GA$-action on the infinite product 
$(X,\mu)=(G/\GA,m_{G/\GA})^{\bZ}$, the index
$[\Rout\Rel_{(X,\GA)}:\OA^*(X,\GA)]$ is infinite countable
\[
  \begin{array}{rllll}
     \OA(X,\GA)&\iso&\Aaut(X,\GA) \ \iso\  S_{\infty}\\
     \OA^*(X,\GA) &\iso& S_{\infty}\times \Gout\GA \\
     \Rout \Rel_{(X,\,\GA)} &\iso& S_{\infty+1}\times
     \Gout\GA
  \end{array}
\]
where $S_\infty$ denotes the full permutation group on
$\bZ$, and $S_{\infty+1}$ the permutation group of $\bZ\cup\{pt\}$
to suggest that the embedding $\OA^*(X,\GA)\subset \Rout\Rel_{X,\GA}$
corresponds to the natural embedding $S_{\infty}\subset S_{\infty+1}$ 
direct product with $\Gout\GA$.
\end{thm}

Let $\Rel$ be an ergodic ${\rm II}_1$-relation on a probability
space $(X,\mu)$, and $E\subset X$ be a measurable subset with
$\mu(E)>0$. The restriction $\Rel_E:=\Rel\cap (E\times E)$ of
$\Rel$ to $E$ is a ${\rm II}_1$-ergodic relation with respect to
the normalized measure $\mu_E:=\mu(E)^{-1}\cdot\mu|_E$. Since
$\Rinn\Rel$ acts transitively on subsets of the same size
(Lemma~\ref{L:full}) for any $F\subset X$ with $\mu(F)=\mu(E)$ the
relation $\Rel_F$ on $(F,\mu_F)$ is isomorphic to $\Rel_E$ on
$(E,\mu_E)$. Hence given a ${\rm II}_1$-relation $\Rel$, for every
$0<t\le 1$ there is a well defined, up to isomorphism, ergodic
${\rm II}_1$-relation $\Rel_t$ obtained from $\Rel$ by restriction
to a subset of measure $t$. (Realizing $(X,\mu)$ as the unit
interval $[0,1]$ one may think of $\Rel_t$ as the restriction of
$\Rel$ to the sub-interval $[0,t]$).

If $\Rel$ has an additional property that $\Rel_t\not\iso\Rel_s$
for all $0<t\neq s \le 1$, one says that $\Rel$ has a
\emph{trivial fundamental group}. Gefter and Golodec \cite{GeG}
proved that 
orbit relations $\Rel=\Rel_{X,\GA}$
generated by free, ergodic, irreducible m.p. actions of higher
rank lattices $\GA$ always have trivial fundamental groups.
(Recent work \cite{Ga} of Gaboriau gives other classes of such
relations).

Regardless whether the fundamental group of $\Rel$ is trivial
or not, all restricted relations $\Rel_t$ obtained from a given
ergodic ${\rm II}_1$-relation $\Rel$ have the same outer
automorphism group: $\Rout\Rel_t\iso\Rout\Rel$ (see
Lemma~\ref{L:induction}). Hence
\begin{thm}\label{T:R-t}
Let $\GA\subset G$ and $(X,\mu,\GA)$ be as in
Theorem~\ref{T:mainN}. For $0<t\le1$ let $\Rel_t$ denote the
(isomorphism class of) equivalence relation obtained from
$\Rel:=\Rel_{X,\GA}$ by a restriction to a subset $E_t\subset X$
of measure $\mu(E_t)=t$. Then $\{\Rel_t\}_{0<t\le1}$ is a family
of mutually non-isomorphic ergodic equivalence relations of type
${\rm II}_1$ with the same outer automorphism group $\Rout\Rel_t
\iso \Rout\Rel_{X,\GA}$. In particular, there exist uncountably
many mutually non-isomorphic ergodic relations with trivial outer
automorphism groups.
\end{thm}
\begin{rems}
\begin{itemize}
\item[(a)]
In \cite{Fu-OE} Thm~D\,(1)--(2) it is shown that for an ergodic
action $(X,\mu,\GA)$ of a lattice $\GA$ in a \emph{simple} higher
rank Lie group $G$, there is a countable set $M_{X,\GA}\subset\bR$
so that for $t\in (0,1)\setminus M_{X,\GA}$ the relation $\Rel_t$
\emph{cannot} be generated by a \emph{free} (mod 0) action of
\emph{any group}. Therefore Theorem~\ref{T:R-t} provides a variety
of examples of such exotic relations without outer automorphisms.
\item[(b)]
In a recent work \cite{MS} Monod and Shalom develop a new
type of "higher rank" superrigidity theorems for products of hyperbolic-like groups.
Using this new tool and the methods of the current paper Monod and Shalom
construct uncountably many \emph{non weakly equivalent} 
relations $\Rel$ of type ${\rm II}_1$ with trivial $\Rout\Rel$ (see \cite{MS} Thm 1.12).
\end{itemize}
\end{rems}


\subsection*{Organization of the paper}
Section~\ref{S:generalities} contains some general facts about
${\rm II}_1$-relations. In section~\ref{S:ME} we discuss the
Measure Equivalence point of view which provides a convenient
framework for the study of $\Rout\Rel_{X,\GA}/\OA^*(X,\GA)$.
Special features of higher rank lattices, especially superrigidity
for cocycles, are used in section~\ref{S:standard-quotients} in a
construction of $\GA$-equivariant \emph{standard quotients}
$\pi:(X,\mu)\to (G/\GA,m_{G/\GA})$ associated to every
$[T]\in\Rout\Rel_{X,\GA}\setminus\OA^*(X,\GA)$, which provide the
proof of Theorem~\ref{T:main}. In section~\ref{S:Ratner} we recall
some ergodic-theoretic applications of Ratner's theorem for
actions on homogeneous spaces. These results are used in
section~\ref{S:main} to assemble the standard quotients for the
proof of Theorem~\ref{T:mainN}, and in sections~\ref{S:examples}
and \ref{S:with-G-GA-factor} to compute the outer automorphism
groups for the standard examples. 
Section~\ref{S:SLnZp} contains the proof of Theorem~\ref{T:SLnZp}.

\section{Generalities}\label{S:generalities}
Let $\Rel$ be an ergodic ${\rm II}_1$ relation on a non-atomic
probability space $(X,\mu)$.
For readers convenience we include the
proof of the following standard fact
\begin{lem}\label{L:full}
For every measurable $E, F\subs X$ with $\mu(E)=\mu(F)>0$
there exists $T\in\Rinn\Rel$ so that $\mu(T E\symdif F)=0$.
\end{lem}
\Pf By \cite{FM} Thm 1, there exists an action $(X,\mu,\GA)$ of
some countable group $\GA$ so that $\Rel=\Rel_{X,\GA}$. Such an
action is necessarily measure-preserving and ergodic. For any
measurable subsets $A,B\subs X$ let $c(A,B):=\sup_\ga\mu(\ga A\cap
B)$. Ergodicity implies that $c(A,B)>0$ whenever $\mu(A)>0$ and
$\mu(B)>0$. Let $E_0:=E$, $F_0:=F$ and define by induction on
$n\ge 1$ measurable sets $E_n, F_n\subs X$ and elements
$\ga_n\in\GA$ as follows: given $E_n, F_n$ choose $\ga_{n}$ so
that
\[
  \mu(\ga_{n}E_n\cap F_n)\ge c(E_n,F_n)/2
\]
and let $E_{n+1}:=E_n\setminus \ga_{n}^{-1}F_n$,
$F_{n+1}:=F_n\setminus \ga_{n}E_n$. Set $E_\infty:=\cap E_n$,
$F_\infty:=\cap F_n$. We have $\mu(E_\infty)=\mu(F_\infty)$
because $\mu(E_n)=\mu(F_n)$ for all finite $n$. In fact
$\mu(E_\infty)=\mu(F_\infty)=0$. Indeed, otherwise one would have
$c(E_n,F_n)\ge c:=c(E_\infty,F_\infty)>0$ for all $n$, contrary to
the choice of $\ga_{n}$ at the stage where $\mu(E_n\setminus
E_{n+1})<c/2$. Hence $E^\prime_n:=E_n\setminus E_{n+1}$ and
$F^\prime_n:=F_n\setminus F_{n+1}$ constitute measurable
partitions of $E$ and $F$ respectively. Defining $T(x)$ to be
$\ga_{n}\cdot x$ if $x\in E^\prime_n$ and $T(x)=x$ for $x\not\in
E$, we get the desired $T\in\Rinn\Rel$. \qed

Given an ergodic ${\rm II}_1$-relation $\Rel$ on $(X,\mu)$, and a
positive measure subset $E\subs X$ we denote by $\Rel_E$ the
restricted relation $\Rel\cap (E\times E)$ on $(E,\mu_E)$, where
$\mu_E=\mu(E)^{-1}\cdot \mu|_E$.
\begin{lem}\label{L:induction}
For a measurable set $E\subs X$ with $\mu(E)>0$
\[
    \Rout\Rel_E\iso\Rout\Rel.
\]
\end{lem}
\Pf
First observe that any $T\in\Raut\Rel_E$ can be extended to a
$\bar{T}\in\Raut\Rel$. To see this choose some measurable partition
$X=E\cup X_1\cup\cdots X_N$ so that $0<\mu(X_i)\le \mu(E)$; and
choose measurable subsets $E_i\subs E$ with $\mu(E_i)=\mu(X_i)$.
By Lemma~\ref{L:full} there exist $S_i, R_i\in\Rinn\Rel$
so that $S_i(X_i)=E_i$ and $R_i(X_i)=T(E_i)$. Define
$\bar{T}$ by $\bar{T}(x)=R_i^{-1}\circ T\circ S_i(x)$ for $x\in X_i$
and $\bar{T}(x)=T(x)$ for $x\in E$ to get a desired $\bar{T}\in\Raut\Rel$.

This extension procedure is well defined on the level of
\emph{outer} classes. In other words if
$\bar{T},\bar{S}\in\Raut\Rel$ are some extensions of some
$T,S\in\Raut\Rel_E$, then $[T]=[S]\in\Rout\Rel_E$ iff
$[\bar{T}]=[\bar{S}]\in\Rout\Rel$. Indeed for $\mu$-a.e. $x\in X$
choose $y\in E$ so that $x\sim y$ and observe that
\[
   \bar{T}(x)\sim \bar{T}(y)=T(y)\rmand S(y)=\bar{S}(y)\sim \bar{S}(x)
\]
Hence $\bar{T}(x)\sim \bar{S}(x)$ for $\mu$-a.e. $x\in X$
iff $T(y)\sim S(y)$ for $\mu_E$-a.e. $y\in E$.

Thus there is a well defined injective map $\Rout\Rel_E\to\Rout\Rel$,
which is easily seen to be a homomorphism of groups. To verify its
surjectivity, note that given any $\bar{T}\in\Raut\Rel$
there is an $\bar{S}\in\Rinn\Rel$ with $\bar{S}(\bar{T}(E))=E$.
Thus $\bar{T}^\prime:=\bar{S}\circ\bar{T}$ maps $E$ to itself,
and $[\bar{T}]=[\bar{T}^\prime]\in\Rout\Rel$
appears as an extension of $[T^\prime|_E]\in\Rout\Rel_E$.
\qed

For the rest of the section we consider a free (mod 0)
ergodic m.p. action $(X,\mu,\GA)$ of some countable group $\GA$,
denoting by $\Rel_{X,\GA}$ the corresponding orbit relation.
\begin{lem}[Gefter \cite{Ge2} Lemmas 2.6, 3.2]\label{L:ICC}
Let $(X,\mu,\GA)$ be a free m.p. ergodic action of a countable
group $\GA$.
\begin{itemize}
\item[{\rm (a)}]
 If $\GA$ has Infinite Conjugacy Classes then
 $\Aaut(X,\GA)\overto{\ep}\OA(X,\GA)$ is an isomorphism.
\item[{\rm (b)}]
 If $\GA$ has the property that any $\ta\in\Gaut\GA$
 with $\ga^\ta=\ga$ on a finite index subgroup $\ga\in\GA_0\subs\GA$
 has to be the identity, then
 \[
    \Ker(\Aaut^*(X,\GA)\overto{\ep}\OA^*(X,\GA))
    =\{ x\mapsto\ga\cdot x\}_{\ga\in\GA}\iso\GA
 \]
\end{itemize}
In particular, the conclusions of (a) and (b) hold for any free ergodic action
$(X,\mu,\GA)$ of an irreducible lattice $\GA$ in a semi-simple Lie
group $G\not\simeq\SL_2(\bR)$.
\end{lem}
\Pf
(a) Any $T\in\Aaut(X,\GA)\,\cap\,\Rinn\Rel_{X,\GA}$ has the form
$T:x\to\xi_{x}\cdot x$ for some measurable $x\mapsto \xi_x\in\GA$
and satisfies $T(\ga\cdot x)=\ga\cdot T(x)$. Hence
\[
    \ga\,\xi_{x}\cdot x=\ga\cdot T(x)=T(\ga\cdot x)
     =\xi_{\ga\cdot x}\,\ga\cdot x
\]
which gives $\xi_{\ga\cdot x}=\ga\,\xi_x\,\ga^{-1}$ because the action is assumed
to be free (mod 0). Thus the distribution $\xi_*\mu$ of $\xi_x$ on $\GA$ is
conjugation invariant, and therefore is uniform
on \emph{finite conjugacy classes} of $\GA$, i.e. supported on $e$. Hence
$T(x)=x$ and $\Ker(\Aaut(X,\GA)\overto{\ep}\OA(X,\GA))$ is trivial.

(b) Any $T\in\Aaut^\ta(X,\GA)\,\cap\,\Rinn\Rel_{X,\GA}$ satisfies
\[
    T(x)=\xi_x\cdot x,\qquad T(\ga\cdot x)=\ga^\ta\cdot T(x)
\]
which gives $\xi_{\ga\cdot x}=\ga^\ta\,\xi_x\,\ga^{-1}$.
For $\xi\in\GA$ let $E_\xi:=\{x\in X\mid \xi_{x}=\xi\}$.
Then $\ga E_\xi=E_{\ga^\ta\xi\ga^{-1}}$.
Observe that for $\xi\neq\xi^\prime\in\GA$ one has
$\mu(E_\xi\cap E_{\xi^\prime})=0$ because the action is free (mod 0).
Hence choosing $\xi_0\in\GA$ with $\mu(E_{\xi_0})>0$
we have $\ga^\ta\xi_0\ga^{-1}=\xi_0$ (equivalently $\xi_0^{-1}\ga^\ta\xi_0=\ga$)
for all $\ga$ in a \emph{finite index} subgroup $\GA_0\subs\GA$. It follows from
the assumption that $\ga^\ta=\xi_0\ga\xi_0^{-1}$ for all $\ga\in\GA$, so that
$T:x\mapsto\xi_0\cdot x$.

Finally, for an irreducible lattice $\GA\subset G\not\simeq\SL_2(\bR)$
the ICC is a standard fact (easy for the group $G$ itself and
follows for $\GA$ using Borel's density theorem), while the condition
for (b) follows from the Strong Rigidity Theorem.
\qed

Given $T\in\Raut\Rel_{X,\GA}$ define a measurable map
$\al_T:\GA\times X\to\GA$ by
\begin{equation}\label{e:aut-cocycle}
   T(\ga\cdot x)=\al_T(\ga,x)\cdot T(x)
\end{equation}
Note that $\al_T(\ga,x)$ is well defined (mod 0) due to
the freeness assumption on the action. Furthermore, one
easily verifies the cocycle property
\[
   \al_T(\ga_2\ga_1,x)=\al_T(\ga_2,\ga_1\cdot x)\,\al_T(\ga_1,x)
\]
for all $\ga_1,\ga_2\in\GA$ and $\mu$-a.e. $x\in X$.
The cocycle $\al_T:\GA\times X\to \GA$ will be called the
\emph{rearrangement cocycle} associated to $T\in\Raut\Rel_{X,\GA}$.
Rearrangement cocycles (as opposed to general ones) have the
following special property:
for $\mu$-a.e. $x\in X$ the correspondence
$\ga\in\GA\ \mapsto\ \al_T(\ga,x)\in\GA$
is a permutation of $\GA$ elements.

Two (general) cocycles $\al,\be:\GA\times X\to\GA$
are said to be \emph{cohomologous in} $\GA$ if there exists
a measurable map $x\mapsto\xi_x\in\GA$, such that
\[
    \al(\ga,x)=\xi_{\ga\cdot x}^{-1}\,\be(\ga,x)\,\xi_{x}
\]
for all $\ga\in\GA$ and $\mu$-a.e. $x\in X$. We denote by
$[\al]_\GA$ the equivalence class of all measurable cocycles
cohomologous (in $\GA$) to $\al$. Note the very special cocycle
$c_1:\GA\times X\to \GA$ given by $c_1(\ga,x)=\ga$, and for a
general $\ta\in\Gaut\GA\ $ let $c_\ta:\GA\times X\to\GA$ stand for
the cocycle $c_\ta(\ga,x)=\ga^\ta$.
\begin{prop}\label{P:basic}
Let $T,S\in\Raut\Rel_{X,\GA}$ be relation automorphisms,
$[T], [S]\in\Rout\Rel_{X,\GA}$ the corresponding classes, and
$\al_T,\al_S:\GA\times X\to\GA$ denote the associated rearrangement
cocycles. Then
\begin{itemize}
\item[{\rm (a)}] $\al_{T\circ S}(\ga,x)=\al_T(\al_S(\ga,x),\,S(x))$.
\item[{\rm (b)}] $\al_T=c_1\iff T\in\Aaut(X,\GA)$.
\item[{\rm (c)}] $\al_T=c_\ta\iff T\in\Aaut^\ta(X,\GA)$.
\item[{\rm (d)}] $[\al_T]_\GA=[c_\ta]_\GA\iff [T]\in \OA^{[\ta]}(X,\GA)$.
\end{itemize}
\end{prop}
\Pf
For $T,S\in\Raut\Rel_{X,\GA}$ compute
\[
    T\circ S(\ga\cdot x)=T(\al_S(\ga,x)\cdot S(x))
     =\al_T(\al_S(\ga,x),\,S(x))\cdot T(S(x))
\]
This proves (a). Statements (b) and (c) follow from the definitions.

Proof of (d). Any $[T]\in\OA^{[\ta]}(X,\GA)$ can be represented
by $T=A\circ J$ where $A\in\Aaut^\ta(X,\GA)$ and $J\in\Rinn\Rel_\GA$ is given by
$J:x\mapsto \xi_{x}^{-1}\cdot x$.
Then for all $\ga\in\GA$ and $\mu$-a.e. $x\in X$
\bqrn
    T(\ga\cdot x)&=&A(\xi_{\ga\cdot x}^{-1}\,\ga\cdot x)
     =\left(\xi_{\ga\cdot x}^{-1}\,\ga\right)^\ta\,\cdot A(x)\\
     &=&\left(\xi_{\ga\cdot x}^{-1}\right)^\ta\,\ga^\ta\,\xi_{x}^\ta\cdot
                              A(\xi_{x}^{-1}\cdot x)
        =\ze_{\ga\cdot x}^{-1}\,\ga^\ta\,\ze_{x}\cdot T(x)
\eqrn
where $\ze_x=(\xi_x)^\ta\in\GA$. Hence
\begin{equation} \label{e:alio}
   \al_T(\ga,x)=\ze_{\ga\cdot x}^{-1}\,\ga^\ta\,\ze_{x}
\end{equation}
and $[\al_T]_\GA=[c_\ta]_\GA$.

On the other hand, assuming that the rearrangement cocycle $\al_T$
associated with $T\in\Raut\Rel_\GA$ satisfies (\ref{e:alio}) for
some $\ta\in\Gaut\GA$ and a measurable $x\mapsto \ze_x\in\GA$, set
$\xi_{x}=(\ze_{x})^{\ta^{-1}}$ and consider the map $A:X\to X$,
defined by $A(x):=\ze_{x}\cdot T(x)$. We have
\bqrn
    A(\ga\cdot x)&=&\ze_{\ga\cdot x}\cdot T(\ga\cdot x)
          =\ze_{\ga\cdot x}\,\ze_{\ga\cdot x}^{-1}\ga^\ta\ze_{x}\cdot T(x)\\
          &=& \ga^\ta\cdot (\ze_{x}\cdot T(x))=\ga^\ta\cdot A(x)
\eqrn
The pushforward measure $A_*\mu$ is absolutely continuous
with respect to $\mu$ (recall that $\GA$ is countable)
and $\GA$-invariant. Ergodicity of the action implies that $A_*\mu=\mu$,
so that $A$ is invertible. Thus $A\in\Aaut^\ta(X,\GA)$, while
the map $J:=A^{-1}\circ T$ is a measure space automorphism. Since
\[
    \xi_{x}\cdot J(x)=\xi_{x}\cdot A^{-1}(T(x))=A^{-1}(\ze_{x}\cdot T(x))=x
\]
the map $J(x)=\xi_{x}^{-1}\cdot x$ is an inner automorphism.
\qed

\section{Measure Equivalence point of view}\label{S:ME}

The following notion of Measure Equivalence Coupling, introduced by
Gromov in \cite{Gr} 0.5.E and considered in \cite{Fu-ME} and
\cite{Fu-OE} by the author, provides a very convenient point of view
on orbit relation automorphisms.
\begin{defn}
A \emph{Measure Equivalence Coupling} of two (infinite) countable groups $\GA$
and $\LA$ is an (infinite) Lebesgue measure space $(\OM,m)$ with two commuting,
free, measure preserving actions of $\GA$ and $\LA$ ,
such that each of the actions has a finite measure fundamental domain.
\end{defn}
We shall use left and right notations for the $\GA$ and $\LA$ actions
\[
    \ga:\om\mapsto \ga\om,\qquad \la:\om\mapsto\om\la
\]
in order to emphasize that the actions commute.
For our current purposes we shall only need \emph{self ME-couplings}
$(\OM,m)$ of $\GA$, i.e. Measure Equivalence Couplings of $\GA$ with itself.
Given such a coupling $(\OM,m)$ let $X, Y\subset \OM$ be
some fundamental domains for the right and the left $\GA$-actions on $(\OM,m)$
respectively. Define the associated measurable maps
\[
     \la=\la_{X}:\GA\times X\to\GA,\qquad
     \rh=\rh_{Y}:Y\times\GA\to\GA
\]
by requiring that for a.e. $x\in X$ (resp. $y\in Y$) one has $\ga
x\in X\la(\ga,x)$ (resp. $y\ga\in\rh(y,\ga) Y$). The left
$\GA$-action on $\OM/\GA$ (resp. the right $\GA$-action on
$\GA\bs\OM$), always denoted by a dot $''\cdot''$, can be identified with
the measure preserving $\GA$-action on $X$ with the finite Lebesgue
measure $m_X=m|_{X}$ (resp. on $Y$ with $m_Y=m|_{X}$) defined by
\[
    \ga\cdot x=\ga\,x\,\la(\ga,x)^{-1},\qquad
    y\cdot\ga=\rh(y,\ga)^{-1}y\ga
\]
With respect to these left and right $\GA$-actions $\la_X$ and $\rh_Y$
become measurable left and right cocycles respectively, namely satisfy:
\[
   \la(\ga_1\ga_2,x)=\la(\ga_1,\ga_2\cdot x)\la(\ga_2,x),\qquad
   \rh(y,\ga_1\ga_2)=\rh(y,\ga_1)\rh(y\cdot\ga_1,\ga_2)
\]
We shall say that a self ME-coupling $(\OM,m)$ is \emph{ergodic}
if the $\GA$-action on $(X,m|_X)$ is ergodic, which is equivalent
to the ergodicity of the $\GA\times\GA$-action on the infinite
space $(\OM,m)$ (see \cite{Fu-ME} Lem~2.2).

With the fundamental domain $X\subset\OM$ for $\OM/\GA$ being fixed,
all fundamental domains $X^\prime\subset\OM$ for $\OM/\GA$
are in one-to-one correspondence with measurable maps $x\mapsto\xi_x\in\GA$:
given a fundamental domain $X^\prime$ one sets $\xi_{x}=\ga$,
if $x\ga\in X^\prime$, while given a measurable $x\mapsto\xi_x\in\GA$
one takes
\[
     X^\prime:=\{x\xi_{x} \mid x\in X\}
\]
The left $\GA$-actions on $X^\prime$ and $X$ are naturally identified
via $\theta:X\to X^\prime$, $\theta : x\mapsto x\xi_{x}$, and the cocycles
$\la_X:\GA\times X\to\GA$, $\la_{X^\prime}:\GA\times X^\prime\to\GA$
are conjugate
\begin{equation}\label{e:change-of-fd}
    \la_{X^\prime}(\ga,\theta(x))
    =\xi_{\ga\cdot x}^{-1}\,\la_X(\ga,x)\,\xi_{x}
\end{equation}
Similar statements hold for the right actions, their fundamental
domains and the associated cocycles.

If $X\subset\OM$ is a fundamental domain for both left and right
$\GA$-actions, we shall say that $X$ is a \emph{two-sided}
fundamental domain.
\begin{lem}[see \cite{Fu-OE} Thm 3.3]\label{L:left-right-fd}
Let $(\OM,m)$ be an ergodic self ME-coupling of some group $\GA$,
and let $X, Y\subset \OM$ be right and left fundamental domains
for $\OM/\GA$ and $\GA\bs\OM$ respectively. Then $\OM$ admits
a two-sided fundamental domain $Z$ iff $m(X)=m(Y)$.
\end{lem}
\Pf Obviously all left fundamental domains have the same
$m$-measure and the same holds for right fundamental domains. Thus
the existence of a two-sided fundamental domain $Z$ implies
$m(X)=m(Z)=m(Y)$. Now assume that $m(X)=m(Y)$. It is well known
that ergodic m.p. actions on finite or infinite Lebesgue spaces
the full group acts transitively on sets of the same measure
(Lemma~\ref{L:full} for the finite measure case). Using the
ergodicity of the $\GA\times\GA$-action on $(\OM,m)$ the condition
$m(X)=m(Y)$ implies that there exist measurable partitions
$X=\bigcup_{i,j} X_{i,j}$, $Y=\bigcup_{i,j} Y_{i,j}$, and elements
$\ga_i^\prime\in\GA$ and $\ga_j^\pprime\in\GA$, so that
$Y_{i,j}={\ga_i^\prime}^{-1} X_{i,j}\ga_j^\pprime$. Then
\[
    \bigcup_{i,j}\ X_{i,j}\ga_j^\pprime\rmand \bigcup_{i,j}\ {\ga_i^\prime}Y_{i,j}
\]
give the same set $Z\subset\OM$. Being formed by piecewise right $\GA$-translates
of $X\iso\OM/\GA$, the set $Z$ is a right fundamental domain for $\OM/\GA$; and
at the same time being formed by piecewise left $\GA$-translates
of $Y\iso\GA\bs\OM$, the same set $Z$ is a left fundamental domain for $\GA\bs\OM$.
\qed

Now consider a free m.p. action $(X,\mu,\GA)$ of some
countable group $\GA$ and let $\Rel_{X,\GA}$ be the corresponding orbit relation.
Given $T\in\Raut\Rel_{X,\GA}$ consider the infinite measure space
$(\OM,m):=(X\times\GA,\mu\times m_\GA)$ with two commuting
$\GA$-actions, as usual written from the left and from the right:
\[
    \ga_1\,(x,\ga):=(\ga_1\cdot x,\al_T(\ga_1,x)\,\ga),\qquad
    (x,\ga)\,\ga_2:=(x,\ga\ga_2)
\]
where $\al_T:\GA\times X\to\GA$ is the rearrangement cocycle
associated with $T\in\Raut\Rel_{X,\GA}$.
The space $(\OM,m)$ with thus defined
$\GA\times\GA$-actions forms an ergodic self ME-coupling of $\GA$,
because $\bar{X}:=X\times\{e_\GA\}\subset\OM$ is a two-sided fundamental
domain. The fact that $\bar{X}$ is a right fundamental domain is obvious.
To see that $\bar{X}$ is a left fundamental domain recall that for a.e. $x\in X$ the map
$\ga\mapsto \al_T(\ga,x)$ is a bijection of $\GA$, so
for $m$-a.e. $(x,\ga_1)$ there is a unique
$\ga\in\GA$ with $\al_T(\ga,x)=\ga_1^{-1}$ which gives
\[
    \ga\,(x,\ga_1)=(\ga\cdot x,\al_T(\ga,x)\ga_1)\in \bar{X}=X\times\{e\}
\]
Also observe that
\begin{equation}\label{e:la-al}
    \la_{\bar{X}}(\ga,x)=\al_T(\ga,x)
\end{equation}
\begin{lem}\label{L:two-sided}
Let $(\OM,m)=(X\times\GA,\mu\times m_\GA)$ be a self ME-coupling
corresponding to $T\in\Raut\Rel_{X,\GA}$.
There is a one-to-one correspondence between two-sided fundamental
domains $\bar{X}^\prime\subset \OM$ and
\begin{equation}\label{e:T-prime}
      T^\prime\in\Raut\Rel_{X,\GA}\rmwith [T^\prime]=[T]\in\Rout\Rel_{X,\GA}
\end{equation}
where $\bar{X}^\prime=\{ (x,\xi_{x}) \mid x\in X\}$ corresponds to
$T^\prime:x\mapsto \xi_{x}^{-1}\cdot T(x)$. Moreover
\[
   \al_{T^\prime}(\ga,x)=\la_{\bar{X}^\prime}(\ga,(x,\xi_{x}))
    =\xi_{\ga\cdot x}^{-1}\la_{\bar{X}}(\ga,(x,e))\,\xi_{x}
    =\xi_{\ga\cdot x}^{-1}\al_T(\ga,x)\,\xi_{x}
\]
\end{lem}
\Pf Suppose that $\bar{X}^\prime\subset\OM=X\times\GA$ is a
two-sided fundamental domain. The fact that both
$\bar{X}=X\times\{e\}$ and $\bar{X}^\prime$ are right fundamental
domains implies that $\bar{X}^\prime$ is of the form
$\{(x,\xi_{x}) \mid x\in X\}$ for some measurable $\xi:X\to \GA$.
In order to verify (\ref{e:T-prime}) for the map $T^\prime:X\to
X$, $T^\prime:x\mapsto \xi_{x}^{-1}\cdot T(x)$, it suffices to
check that $T^\prime$ is one-to-one (mod 0), the relations between
the cocycles $\al_{T^\prime}$, $\la_{\bar{X}^\prime}$,
$\la_{\bar{X}}$ and $\al_T$ being straightforward.

Assume that $T^\prime(x)=T^\prime(y)$ which means
$\xi_{x}^{-1}\cdot T(x)=\xi_y^{-1}\cdot T(y)$.
Then $T(x)$ and $T(y)$ are on the same $\GA$-orbit in $X$,
and so are $x$ and $y$,
i.e. $y=\ga\cdot x$ for some $\ga\in\GA$. Thus
\[
   \xi_{x}^{-1}\cdot T(x)=\xi_{\ga\cdot x}^{-1}\cdot T(\ga\cdot x)
   =\xi_{\ga\cdot x}^{-1}\,\al_T(\ga,x)\cdot T(x)
\]
which means that $\xi_{\ga\cdot x}=\al_T(\ga,x)\,\xi_{x}$. In $\OM$ we have
\[
    \ga (x,\xi_{x})=(\ga\cdot x,\al_T(\ga,x)\xi_{x})
    =(\ga\cdot x,\xi_{\ga\cdot x})
\]
with both $(x,\xi_{x})$ and $(\ga\cdot x,\xi_{\ga\cdot x})$
in $\bar{X}^\prime$.
Since $\bar{X}^\prime$ is a two-sided fundamental domain, in particular a
left fundamental domain,
it follows that $\ga=e$ and $x=y$.
Hence $T^\prime$ is indeed a measure space automorphism of $(X,\mu)$ and
the rest of its properties follow automatically.
The fact that $T^\prime$ as in (\ref{e:T-prime}) gives rise to a two-sided
fundamental domain $\bar{X}^\prime$ is proved by back tracking the above argument.
\qed

Next consider an equivariant quotient map $\PH:(\OM,m)\to(\OM_0,m_0)$
of self ME-couplings of $\GA$, i.e. a measurable map $\PH:\OM\to\OM_0$ such that
\[
   \PH_*m=m_0\rmand \PH(\ga_1\om\ga_2)=\ga_1\,\PH(\om)\,\ga_2
\]
Observe that the preimage $X:=\PH^{-1}(X_0)$ (resp. $Y:=\PH^{-1}(Y_0)$)
of any right fundamental domain $X_0\subset\OM_0$ (resp. any left fundamental domain $Y_0\subset\OM_0$)
is a right (resp. left) fundamental domain in $\OM$. If $X=\PH^{-1}(X_0)$
we shall say that $X\subset\OM$ and $X_0\subset\OM_0$ are
$\PH$-\emph{compatible}. Note also that if $(\OM,m)$ is an
ergodic coupling then so is $(\OM_0,m_0)$, and if $(\OM,m)$ admits a
two-sided fundamental domain then
\[
   m_0(X_0)=m(X)=m(Y)=m_0(Y_0)
\]
so that $(\OM_0,m_0)$ also admits a two-sided fundamental domain $Z_0$,
and taking $Z:=\PH^{-1}(Z_0)$ we obtain a \emph{two-sided fundamental domain}
for $(\OM,m)$ which is $\PH$-compatible with $Z_0\subset\OM_0$.

Observe that for $\PH$-compatible right fundamental domains $X\subset\OM$
and $X_0\subset\OM_0$ one has
\[
   \la_{X}(\ga,\om)=\la_{X_0}(\ga,\PH(\om))
\]
Realizing the natural left $\GA$-action on $(\OM,m)/\GA$ by the
$\GA$-action
\[
   \ga: x \mapsto \ga\cdot x=\ga\, x\,\la_{X}(\ga,x)^{-1}
\]
on a $\PH$-compatible fundamental domain $X\subset\OM$,
one obtains a $\GA$-equivariant quotient map $X\overto{\PH} X_0$
which is a concrete realization of the left $\GA$-equivariant
map $(\OM,m)/\GA\to (\OM_0,m_0)/\GA$ defined by $\PH$.
This discussion is summeraized by the following
\begin{prop}\label{P:ME}
Let $(X,\mu,\GA)$ be a free, ergodic, measure preserving action,
$T\in\Raut\Rel_{X,\GA}$ and let $(\OM_{[T]},m)$ be the
corresponding self ME-coupling of $\GA$. Assume that $(\OM_{[T]},m)$
has an equivariant quotient ME-coupling
$\PH:(\OM_{[T]},m)\to (\OM_0,m_0)$.
Fix a two sided fundamental domain $X_0\subset\OM_0$,
denote by $(X_0,\mu_0,\GA)$ the left $\GA$-action on
$(X_0,\mu_0)\iso(\OM_0,m_0)/\GA$, and let
\[
    \pi:(X,\mu,\GA)\to (X_0,\mu_0,\GA)
\]
denote the $\GA$-equivariant quotient map induced by $\PH$. Then
there exists a
$\hat{T}\in\Raut\Rel_{X,\GA}$ with $[T]=[\hat{T}]\in\Rout\Rel_{X,\GA}$
so that
\[
     \al_{\hat{T}}(\ga,x)=\la_{X_0}(\ga,\pi(x)).
\]
\end{prop}

\section{Superrigidity and Standard Quotients} \label{S:standard-quotients}
In this section we specialize to actions of irreducible lattices $\GA$
in higher rank semi-simple Lie groups $G$.
\begin{prop}[see \cite{Fu-ME} Thm 4.1]\label{P:ME-quotients}
Let $G$ be a semi-simple, connected, center-free, real Lie group
without non-trivial compact factors and with $\rk(G)\ge 2$. 
Let $\GA\subset G$ be an irreducible lattice and $(X,\mu,\GA)$ be a
measure preserving, ergodic, irreducible, 
essentially free $\GA$-action. 
Given any $T\in\Raut\Rel_{X,\GA}$ let
$(\OM_{[T]},m)$ be the associated self ME-coupling as in
section~\ref{S:ME}. Then there exists a well defined class
$[\ta]\in\Gout G $ so that given any representative $\ta$ of
$[\ta]$ there exists a measurable map $\PH:\OM_{[T]}\to G$ defined
uniquely (mod 0) so that
\[
   \PH(\ga_1\om\ga_2)=\ga_1^\ta\PH(\om)\ga_2\qquad(\ga_1,\ga_2\in\GA)
\]
and one of the following two alternatives holds:
\begin{itemize}
\item[{\rm (a)}]
  either $\PH_*m$ coincides with the Haar measure $m_G$ on $G$,
  normalized so that $\GA$ has covolume one, or
\item[{\rm (b)}]
  $\PH_*m$ is an atomic measure of the form
  \[
  	\frac{1}{k} \sum_{i=1}^k \sum_{\ga\in\GA} \de_{g_i\ga}
  \] 
  where $\{g_i\}_1^k\subset G$ are such that $\{g_1\GA,\dots,g_k\GA\}$
  is a single finite $\ta(\GA)$-orbit on $G/\GA$. 
  In particular, $\GA$ has a subgroup $\GA_1$ of index $k$ so that
  $\ta(\GA_1)$ has index $k$ in $g_1\GA g_1^{-1}$, and $\ta(\GA)$ and $\GA$ 
  are commensurable.
\end{itemize}
If the $\GA$-action on $(X,\mu)$ is aperiodic, then either {\rm (a)} holds
or in alternative {\rm (b)} we have $k=1$ which means that
\begin{itemize}
\item[{\rm (b')}]
  $\PH_*m$ coincides with the counting measure $m_{\GA^\prime}$
  on $\GA^\prime=\ta(\GA)\subset G$ where  
  $\ta(\GA)=g\GA g^{-1}$ for some $g\in G$.
\end{itemize}
\end{prop}
This proposition is essentially \cite{Fu-ME} Thm 4.1, the proof of
which is based on Zimmer's superrigidity for cocycles and Ratner's
theorem. 
In \cite{Fu-ME} the statement is formulated in a slightly
different form and only for lattices in higher rank \emph{simple}
Lie groups. Since we need some details of the proof to be used
later, we include the main arguments here.

Fix a $T\in\Raut\Rel_{X,\GA}$ representing $[T]$ and consider the
rearrangement cocycle $\al_T:\GA\times X\to\GA\subset G$ as a
$G$-valued cocycle. This cocycle is Zariski dense in $G$ (this is
a form of Borel's density theorem, see \cite{Zi-book} p. 99, or
\cite{Fu-ME} Lemma 4.2). Thus the assumption that $\GA$ is a higher
rank lattice with an \emph{irreducible} action on $(X,\mu)$ allows
to apply Zimmer's superrigidity for measurable cocycles theorem
\cite{Zi-book} (in \cite{Fu-ME} Thm 4.1 we did not use
irreducibility of the action and therefore had to restrict the
discussion to lattices in higher rank \emph{simple} groups $G$).
Hence there exists a Borel map $\ph:X\to G$ and a homomorphism
$\ta:\GA\to G$, so that
\begin{equation}\label{e:al-ph}
    \al_T(\ga,x)=\ph(\ga\cdot x)^{-1}\,\ga^\ta\,\ph(x)
\end{equation}
for $\ga\in\GA$ and $\mu$-a.e. $x\in X$. By Margulis' superrigidity
$\ta$ extends to a $G$-automorphism and we denote by $\ta\in\Gaut G $
this extension. Defining the map $\PH:\OM_{[T]}=X\times\GA\to G$ by
\begin{equation}\label{e:ph-PH}
     \PH(x,\ga):=\ph(x)\ga
\end{equation}
one verifies
\bqrn
    \PH(\ga_1 (x,\ga) \ga_2)&=&\PH(\ga_1\cdot x,\,\al_T(\ga_1,x)\ga\ga_2)
    =\ph(\ga_1\cdot x)\al_T(\ga_1,x)\ga\ga_2\\
    &=&\ph(\ga_1\cdot x)\,\ph(\ga_1\cdot x)^{-1}\ga_1^\ta\,\ph(x)\,\ga\,\ga_2\\
    &=&\ga_1^\ta\,\PH(x,\ga)\,\ga_2
\eqrn

Choose $F\subset G$ a Borel fundamental domain for $G/\GA$ and let
$\bar{X}^\prime:=\PH^{-1}(F)$. 
Hence $\bar{X}^\prime\subset\OM_{[T]}$ is a fundamental domain for
$\OM_{[T]}/\GA$ so that $m(\bar{X}^\prime)=1$. 
This implies that the pushforward measure $m_0:=\PH_* m$ on $G$ 
has $m_0(F)=1$ (in particular $m_0$ is finite on compact sets) while the restriction
$m_0|_F$ defines a regular Borel probability measure $\mu_0$ on $G/\GA$, which is
invariant and ergodic for the left $\ta(\GA)$-action.

An application of Ratner's theorem (see \cite{Fu-ME} Lem~4.6 with
an easy modification needed to handle semi-simple rather simple
Lie groups) implies that $\mu_0$ is either (i) $\mu_0=m_{G/\GA}$ - the normalized 
Haar measure $m_{G/\GA}$, or (ii) is an atomic measure.

In case (i) the map $\PH$ defined in
(\ref{e:ph-PH}) clearly maps $m$ on $\OM_{[T]}$ to the Haar
measure $m_G$ as in Proposition~\ref{P:ME-quotients}\,(a). 
The uniqueness statements in Proposition~\ref{P:ME-quotients}
follow from \cite{Fu-ME} Thm 4.1.

In case (ii) the atomic $\ta(\GA)$-invariant measure $\mu_0$ on
$G/\GA$ has to be concentrated on a single finite $\ta(\GA)$-orbit 
$\{g_1\GA,\dots,g_k\GA\}$ with equal weights $1/k$. 
Let $\GA_1$ be the stabilizer of $g_1\GA\in G/\GA$.
Then $[\GA:\GA_1]=k$ and $\ta(\GA_1)g_1\GA=g_1\GA$ i.e.
$\ta(\GA_1)$ has index $k$ in $g_1\GA g_1^{-1}$.

The preimage $\OM_1=\PH^{-1}(g_1\GA)$ is $\GA_1\times\GA$-invariant 
set which gives rise to a measurable $\GA_1$-invariant subset $X_1$ of $X$
with $\mu(X_1)=1/k$.
If $\GA$-action on $(X,\mu)$ is aperiodic, then necessarily $k=1$
and $m_0=\sum_{\ga\in\GA}\de_{g\ga}$ and $\ta(\GA)=g\GA g^{-1}$.
\qed
\begin{rem}
The uniqueness of $\PH_{[T]}$ in particular implies that
the rearrangement cocycle $\al_T$ can be written in the form (\ref{e:al-ph})
with the measurable map $\ph:X\to G$ being uniquely defined (mod 0)
as soon as a representative $\ta\in\Gaut G $ of $[\ta]\in\Gout G $ is chosen.
Hereafter this unique ``straightening'' map $\ph$ will be denoted by $\ph_{T,\ta}$.
\end{rem}

\begin{thm}[{\bf Standard Quotients}]\label{T:standard-quotients}
Let $G$ be a semi-simple, connected, center-free, real Lie group
without non-trivial compact factors and with $\rk(G)\ge 2$; 
$\GA\subset G$ - an irreducible lattice and $(X,\mu,\GA)$ be a
measure preserving, ergodic, irreducible, essentially free $\GA$-action. 
Then every $[T]\in\Rout\Rel_{X,\GA}$ defines a unique class $[\ta]\in\Gout G $ 
such that given any representative $\ta\in\Gaut G $ of $[\ta]$ there is a measurable
map $\pi:X\to G/\GA$, defined uniquely (mod 0) and satisfying 
\[
	\pi(\ga\cdot x)=\ta(\ga)\cdot\pi(x)
\] 
for $\mu$-a.e. $x\in X$ and all $\ga\in\GA$. There are two alternatives:

\noindent{either the following equivalent conditions hold}:
\begin{itemize}
  \item[{\rm (a1)}] the distribution of $\ph_{T,\ta}(x)$ on $G$
     is absolutely continuous with respect to the Haar measure $m_G$;
  \item[{\rm (a2)}] $\pi_*\mu=m_{G/\GA}$ - the $G$-invariant
     probability measure on $G/\GA$;
  \item[{\rm (a3)}] there exists $\hat{T}\in\Raut\Rel_{X,\GA}$ with
     $[\hat{T}]=[T]$ and $\pi(x)=\ph_{\hat{T},\ta}(x)\,\GA$.
\end{itemize}
or the following equivalent conditions hold:
\begin{itemize}
  \item[{\rm (b1)}] the distribution of $\ph_{T,\ta}(x)$ on $G$
       is purely atomic;
  \item[{\rm (b2)}] $\pi_*\mu=k^{-1}\sum_1^k\de_{g_i\GA}$ 
       where  $\{g_1\GA,\dots,g_k\GA\}$ is a finite $\ta(\GA)$-orbit on $G/\GA$;
       $\GA$ contains a subgroup $\GA_1$ of index $k$ so that
       $\ta(\GA_1)$ is a subgroup of index $k$ in $g_1\GA g_1^{-1}$;
       and $X_1=\pi^{-1}(\{g_1\GA\})$ is a $\GA_1$-ergodic 
       components of $(X,\mu)$ with $\mu(X_1)=1/k$;
   \item[{\rm (b3)}] there exists $\hat{T}\in\Raut\Rel_{X,\GA}$ with
     $[\hat{T}]=[T]$ and 
     \[
     	\ph_{\hat{T},\ta}(x)=g_1 \quad{\rm for\ }\mu-{\rm a.e.\ } x\in X_1\subset X.
     \]
\end{itemize}
If $\GA$-action on $(X,\mu)$ is aperiodic then 
conditions {\rm (a1)-(a3)} above are equivalent to
\begin{itemize}
  \item[{\rm (a4)}] $[T]\not\in\OA^*(X,\GA)$,
\end{itemize}
while their alternatives {\rm (b1)--(b3)} are equivalent to
\begin{itemize}
  \item[{\rm (b4)}] $[T]\in\OA^*(X,\GA)$;
\end{itemize}
moreover in {\rm (b2)--(b3)} one has $k=1$ and these conditions 
take the following form:
\begin{itemize}
  \item[{\rm (b2')}] $\pi_*\mu=\de_{g\GA}$ 
       where  $g\in G$ satisfies $\ta(\GA)=g\GA g^{-1}$;
   \item[{\rm (b3')}] there exists $\hat{T}\in\Raut\Rel_{X,\GA}$ with
     $[\hat{T}]=[T]$ and 
     \[
     	\ph_{\hat{T},\ta}(x)=g\qquad{\rm for\ }\mu{\rm -a.e.}\  x\in X.
     \]
\end{itemize}
\end{thm}

\Pf Consider the self ME-coupling $(\OM_{[T]},m)$ with the
corresponding outer class $[\ta]\in\Gout G $. Given a choice
$\ta\in\Gaut G$ of $[\ta]$ let 
\[
	\PH:\OM_{[T]}\to G
\] 
be the $\ta(\GA)\times\GA$-equivariant map as in
Proposition~\ref{P:ME-quotients}. Then $\PH$ uniquely defines a
measurable map
\[
    \pi:(X,\mu,\GA)\iso(\OM_{[T]},m)/\GA\to G/\GA,\qquad
    \pi(\ga\cdot x)=\ga^\ta\cdot \pi(x)
\]
Let us show that the alternatives (a) and (b) in
Proposition~\ref{P:ME-quotients} yield mutually exclusive
conditions (a1)--(a3) and (b1)--(b3) respectively.

In case (a) where $\PH_*m=m_G$, (a1)--(a3) follow
from Proposition~\ref{P:ME} and the construction
(\ref{e:ph-PH}) of $\PH$.

Case (b): $\PH_*m=k^{-1}\sum_{i=1}^k \sum_{\ga\in\GA} \de_{g_i\ga}$ 
where $\{g_1\GA,\dots,g_k\GA\}$ is a single $\ta(\GA)$-orbit on $G/\GA$.
Condition (b1) is clearly satisfied. 
Let 
\[
	\GA_i=\{\ga\in\GA\,\mid\, \ga^\ta g_i\GA=g_i\GA\}
	\rmand X_i=\pi^{-1}(\{g_i\GA\}).
\]
where $\pi:X\to\{g_1\GA,\dots,g_k\GA\}$ 
is the $\GA$-equivariant map above.
Then conjugate groups $\GA_i$ have index $k$ in $\GA$, and 
$\GA$ permutes the disjoint sets $X_i$ (and so $\mu(X_i)=1/k$)
while each $X_i$ is $\GA_i$-invariant for $i=1,\dots,k$. 
Moreover $\GA_i$ acts ergodically on $X_i$
because $\Rel_{X_i,\GA_i}=\Rel_{X,\GA}\cap (X_i\times X_i)$.
This proves (b2). 

The set $X_0=\{g_1,\dots,g_k\}$ forms a fundamental domain for
the $\ta(\GA)\times\GA$-action on $G$.
The corresponding cocycle $\la_{X_0}$ satisfies
\[
	\la_{X_0}(\ga_1, g_1)=g_1^{-1} \ga_1^\ta g_1\qquad (\ga_1\in\GA_1)
\]
Applying Proposition~\ref{P:ME} we obtain 
$\hat{T}\in\Raut\Rel_{X,\GA}$ with $[\hat{T}]=[T]\in\Rout\Rel_{X,\GA}$ and
\begin{equation}\label{e:g1}
   	\al_{\hat{T}}(\ga,x)=\la_{X_0}(\ga,\pi(x))=g_1^{-1}\,\ga^\ta\, g_1
\end{equation}
for all $\ga\in\GA_1$ and a.e. $x\in X_1=\pi^{-1}(\{g_1\GA\})=\PH^{-1}(\{g_1\})$.
We deduce that $\ph_{\hat{T},\ta}(x)=g_1$ for $x\in X_1$, proving (b3).

If the $\GA$-action $(X,\mu)$ is aperiodic, one has $k=1$ so that
(b2), (b3) take the form of (b2'), (b3'). 
Condition (b3) follows from (\ref{e:g1}) and Proposition~\ref{P:basic}\,(c).
The latter also explains why (b4) is incompatible with (a1)--(a3).
\qed

\noindent{\bf Proof fo Theorem~\ref{T:main}.}
For $\ta\in\Gaut G$ the $\GA$-action on $G/\ta^{-1}(\GA)$ 
is isomorphic 
to the $\ta$-twisted 
$\GA$-action on $G/\GA$, both with the Haar measure. 
Since $(X,\mu,\GA)$ is assumed not to have these actions among its  
measurable quotients, any $T\in\Raut\Rel_{X,\GA}$ fails 
condition ${\rm (a2)}$ in Theorem~\ref{T:standard-quotients}, 
while satisfies the alternatives (b1-4), which means that $[T]\in\OA^*(X,\GA)$.
\qed

\section{Some applications of Ratner's Theorem}\label{S:Ratner}
In this section we recall some applications of Ratner's Theorem
(see \cite{Ra} and references therein). Note that in these results
there are no restrictions on the rank of the semi-simple group $G$.
In fact the results remain true whenever 
$G$ is a connected Lie group generated by $\Ad$-unipotent elements
and $\GA\subset G$ is a closed
subgroup so that $G/\GA$ carries a $G$-invariant probability measure.
\begin{thm}[cf. Ratner \cite{Ra} Thm E2]
  \label{T:quotients}
Let $\GA$ be an irreducible lattice in a semi-simple connected
real Lie group $G$, $\LA$ and $\LA^\prime$ be lattices in some connected
Lie groups $H$ and $H^\prime$, $\rh:G\to H$ and $\rh^\prime:G\to H^\prime$
be continuous homomorphisms such that the $\GA$-actions
\[
   \ga\,:\,h\LA\mapsto\rh(\ga)h\LA,\qquad
   \ga\,:\,h^\prime\LA^\prime\mapsto\rh^\prime(\ga)h^\prime\LA^\prime
\]
on $(H/\LA,m_{H/\LA})$ and $(H^\prime/\LA^\prime,m_{H^\prime/\LA^\prime})$
are ergodic. Assume that there exists a measurable $\GA$-equivariant
quotient map
\[
   \pi:(H/\LA,m_{H/\LA})\overto{} (H^\prime/\LA^\prime,m_{H^\prime/\LA^\prime})
\]
Then there exists a $t\in H^\prime$ and a surjective continuous homomorphism
$\si:H\to H^\prime$ such that
\begin{itemize}
\item[{\rm (i)}]
    $\si(\LA)$ is a finite index subgroup of $\LA^\prime$,
\item[{\rm (ii)}]
    $\pi(h\LA)=t\si(h)\LA^\prime$ for a.e. $h\in H$,
\item[{\rm (iii)}]
    $\si\circ\rh(\ga)=t\rh^\prime(\ga)t^{-1}$ for $\ga\in\GA$.
\end{itemize}
If $\pi$ is one-to-one then $\si:H\to H^\prime$ is an isomorphism
and $\si(\LA)=\LA^\prime$. In particular, for the above
$\GA$-action on $(H/\LA,m_{H/\LA})$
\[
   \Aaut^*(H/\LA,m_{H/\LA},\GA)\iso N_{\Aff(H/\LA)}(\rh(\GA))
\]
\end{thm}
In \cite{Wi} Witte considers a more general question of a classification
of all measurable equivariant quotients $(H/\LA,m_{H/\LA})\to (Y,\nu)$
showing that $(Y,\nu)$ has an algebraic description (slightly more general
than $H^\prime/\LA^\prime$ as above).
However Theorem~\ref{T:quotients} suffices for our purposes.
It is deduced from the more general
Theorem~\ref{T:joinings} below by
considering the measure $\nu$ on $H/\LA\times H^\prime/\LA^\prime$
obtained by the lift of $m_{H/\LA}$ to the graph of
$\pi:H/\LA\to H^\prime/\LA^\prime$.
\begin{thm}[cf. Ratner \cite{Ra} Thm E3]
   \label{T:joinings}
Let $\GA\subset G$, $\LA\subset H$, $\LA^\prime\subset H^\prime$,
$\rh:G\to H$ and $\rh^\prime:G\to H^\prime$ be as in Theorem~\ref{T:quotients}.
Let $\nu$ be a probability measure on $H/\LA\times H^\prime/\LA^\prime$
which projects onto $m_{H/\LA}$ and $m_{H^\prime/\LA^\prime}$, and is invariant
and ergodic for the diagonal $\GA$-action
\[
   \ga\,:\,(h\LA,h^\prime\LA^\prime)\,\mapsto\,
           (\rh(\ga)h\LA,\rh^\prime(\ga)h^\prime\LA^\prime)
\]
Then there exist closed normal subgroups $N\normal H$,
$N^\prime\normal H^\prime$, an element $t\in H^\prime/N^\prime$
and a continuous isomorphism $\si_1$ from $H_1:=H/N$ to
$H_1^\prime:=H^\prime/N^\prime$, so that
\begin{itemize}
\item[{\rm (i)}]
    $\LA_1:=\LA N\subset H_1$ and
    $\LA_1^\prime:=\LA^\prime N^\prime\subset H_1^\prime$
    are lattices,
\item[{\rm (ii)}]
    There are finite index subgroups $\DE_1\subs\LA_1$, $\DE_1^\prime\subs\LA_1^\prime$
    so that $\si_1(\DE_1)=\DE_1^\prime$,
\item[{\rm (iii)}]
    $\si_1\circ\rh(\ga)=t\rh^\prime(\ga)t^{-1}$ for $\ga\in\GA$,
\item[{\rm (iv)}]
    The measure $\nu$ is $N\times N^\prime$-invariant and its projection
    $\nu_1$ to $H_1/\LA_1\times H_1^\prime/\LA_1^\prime$ can be obtained from
    the lift $m_f$ of $m_{H_1/\DE_1}$ to the graph of
    $H_1/\DE_1\overto{f} H_1^\prime/\DE_1^\prime$ where
    $f(h\DE_1)=t\si_1(h)\DE_2$,
    by $\nu_1=p_* m_f$ where $p$ is a finite-to-one projection
    \[
       (H_1/\DE_1)\times (H_1^\prime/\DE_1^\prime)\overto{p}
       (H_1/\LA_1)\times (H_1^\prime/\LA_1^\prime).
    \]
\end{itemize}
\end{thm}
Theorem \cite{Ra}~E3 and its corollary \cite{Ra}~E2 were proved by
M.~Ratner as an application of the main theorem (\cite{Ra} Thm 1).
In all these results the acting group is assumed to be generated
by $\Ad$-unipotent elements. In order to deduce the results for
actions of lattices $\GA\subset G$, needed for our purposes, one
uses the suspension construction replacing the $\GA$-invariant
measure $\nu$ on $H/\LA\times H^\prime/\LA^\prime$ by the
$G$-invariant measure $\tilde{\nu}$ on $G\times_\GA H/\LA\times
H^\prime/\LA^\prime$ and applying Ratner's classification of
invariant measures (\cite{Ra} Thm 1) to the action of the
semi-simple group $G$ which is generated by $\Ad$-unipotents. The
reader is referred to the paper \cite{Sh} of Shah (Corollary~1.4)
or Witte (\cite{Wi} proof of Corollary~5.8) for the precise
argument.

\section{Proofs of Theorem~\ref{T:mainN} and Corollary~\ref{C:entropy}} \label{S:main}

\subsection*{Proof of Theorem~\ref{T:mainN}}

Let $\GA\subset G$ and $(X,\mu,\GA)$ be as in
Theorem~\ref{T:mainN}, and let
$n:=[\Rout\Rel_{X,\GA}:\OA^*(X,\GA)]\in\{1,2,\dots,\infty\}$. If
$n=1$ there is nothing to prove. If $1< n\le\infty$ set $T_0=\Id$
and choose representatives $T_i\in\Raut\Rel_{X,\GA}$, $1\le i<n$,
for the cosets $\OA^*(X,\GA)\bs\Rout\Rel_{X,\GA}$. In other words
choose $T_i$ so that for $0\le i\neq j<n$ we have
\[
   [T_i][T_j]^{-1}\not\in\OA^*(X,\GA)
\]
Since $[T_i]\not\in\OA^*(X,\GA)$ for $1\le i<n$, by
Theorem~\ref{T:standard-quotients}\,(a) there are  $\ta_i\in\Gaut G $
and measurable maps $\pi_i:X\to G/\GA$ satisfying
\[
    (\pi_i)_*\mu=m_{G/\GA},\qquad \pi_i(\ga\cdot x)=\ga^{\ta_i}\cdot\pi_i(x)
\]
It remains to prove that the map
\[
   \pi:X\overto{}\prod_{i=1}^{n-1} G/\GA,\qquad\pi(x)=(\pi_1(x),\pi_2(x),\dots)
\]
takes $\mu$ onto the product measure
$m_{G^{n-1}/\GA^{n-1}}=\prod_{i=1}^{n-1} m_{G/\GA}$.
We shall prove by induction on \emph{finite} $k$ in the range $1\le k<n$
that the map $\pi^{(k)}(x):=(\pi_1(x),\dots,\pi_k(x))$ satisfies
\begin{equation}\label{e:induction}
    \pi^{(k)}_*\mu=m_{G^k/\GA^k}
\end{equation}
(Note that this is sufficient even if $n=\infty$ because the
infinite product measure is determined by its values on finite
cylinder sets). The case $k=1$ is covered by
Theorem~\ref{T:standard-quotients}\,(a2). 
Assuming (\ref{e:induction}) for $k-1$ we apply Theorem~\ref{T:joinings} to
\[
   \begin{array}{lll}
   H:=G^{k-1}\qquad &\LA:=\GA^{k-1}\qquad
      &\rh:=\ta_1\times\cdots\times\ta_{k-1}\\
   H^\prime:=G\qquad &\LA^\prime:=\GA\qquad &\rh^\prime:=\ta_k
   \end{array}
\]
and the probability measure $\nu:=\pi^{(k)}_*\mu$ on $H/\LA\times
H^\prime/\LA^\prime=G^k/\GA^k$. By the induction hypothesis $\nu$
projects onto $m_{H/\LA}$ in the first factor, and as
$[T_k]\not\in\OA^*(X,\GA)$, $\nu$ projects onto
$m_{H^\prime/\LA^\prime}$ in the second factor. If $N=H=G^{k-1}$
then necessarily $N^\prime=H^\prime=G$, so that
\[
   \nu=m_{H/\LA}\times m_{H^\prime/\LA^\prime}=m_{G^k/\GA^k}
\]
proving the induction step.

It remains to show that the other alternative, namely $N\normal
G^{k-1}$ and $N^\prime\normal G$ being \emph{proper} normal
subgroups, cannot occur. By Theorem~\ref{T:joinings}\,(i)
$\LA_1=\GA N^\prime\subset G/N^\prime$ forms a lattice in
$G/N^\prime$ which is possible only if $N^\prime=\{e\}$ because
$\GA\subset G$ is irreducible. Thus $N\normal G^{k-1}$ is such
that $G^{k-1}/N\iso G$ and $\GA^{k-1} N$ forms a lattice in
$G^{k-1}/N\iso G$. This means that for some $j\in\{1,\dots,k-1\}$
\bqrn
   &&N=\{ (g_1,\dots,g_{k-1})\in G^{k-1} \mid g_j=e \}\\
   &&\si_1((g_1,\dots,g_{k-1})N)=\si(g_j)
\eqrn where $\si\in \Gaut G$ is such that for some $t\in G$,
$\si\circ\ta_j(g)=t\ta_k(g) t^{-1}$ and $\si(\DE)=\DE^\prime$ for
some finite index subgroups $\DE,\DE^\prime\subs\GA$. In this case
the distribution $\nu_1$ of the pairs $(\pi_j(x),\pi_k(x))$ on
$G/\GA\times G/\GA$ is a projection under the finite-to-one map
\[
   G/\DE\times G/\DE^\prime\,\overto{}\,G/\GA\times G/\GA
\]
of the measure $m_f$ which is a lift of $m_{G/\DE}$ to the graph of
\[
    f:G/\DE\overto{} G/\DE^\prime,\qquad
    f(g\DE)=t\si(g)\DE^\prime
\]
By Theorem~\ref{T:standard-quotients}\,(a3) there exist
$\hat{T}_j$ and $\hat{T}_k\in\Raut\Rel_{X,\GA}$
with $[\hat{T}_j]=[T_j]$, $[\hat{T}_k]=[T_k]$ so that
for $i=j,k$ the rearrangement cocycles
\[
   \al_i:=\al_{\hat{T}_i}\ :\GA\times X\overto{} \GA
\]
satisfy $\al_i(\ga,x)=\ph_i(\ga\cdot x)^{-1}\,\ga^{\ta_i}\,\ph_i(x)$
with $\pi_i(x)=\ph_i(x)\GA$.
The structure of the distribution $\nu_1$ of $(\pi_j(x),\pi_k(x))$
described above implies that the distribution of
$(\ph_j(x)^{\si})^{-1}\ph_k(x)$ on $G$ is \emph{purely atomic}.
Let $S:=\hat{T}_k\circ \hat{T}_j^{-1}\in\Raut\Rel_{X,\GA}$ and let
$\si^\prime\in\Gaut G$ and $\ps=\ph_{S,\si^\prime}:X\to G$ be such that
\[
     \al_S(\ga,x)=\ps(\ga\cdot x)^{-1}\,\ga^{\si^\prime}\,\ps(x)
\]
Applying Proposition~\ref{P:basic}\,(a) to $\hat{T}_k=S\circ
\hat{T}_j$ we obtain that for all $\ga\in\GA$ and $\mu$-a.e. $x\in
X$ \bqrn
    &&\ph_k(\ga\cdot x)^{-1}\,\ga^{\ta_k}\,\ph_k(x)=
       \al_k(\ga,x)=\al_S\left(\al_j(\ga,x),\hat{T}_j(x)\right)\\
    &&=\ps\left(\al_j(\ga,x)\cdot \hat{T}_j(x)\right)^{-1}
         \al_j(\ga,x)^{\si^\prime}\ \ps(\hat{T}_j(x))\\
    &&=\ps\left(\hat{T}_j(\ga\cdot x)\right)^{-1}\,
       \left(\ph_j(\ga\cdot x)^{-1}\,\ga^{\ta_j}\ \ph_j(x)\right)^{\si^\prime}\,
       \ps(\hat{T}_j(x))\\
    &&=\left(\ph_j(\ga\cdot x)^{\si^\prime} \ps(\hat{T}_j(\ga\cdot x))\right)^{-1}\,
       \ga^{\si^\prime\circ\ta_j}\,\,
       \left(\ph_j(x)^{\si^\prime}\ps(\hat{T}_j(x))\right)
\eqrn
Replacing $\si^\prime$ by $\si\in\Gaut G$ (so that
$\ta_k=\si\circ\ta_j$) and changing $\ps=\ph_{S,\si^\prime}$ to
$\ph_{S,\si}$ accordingly, we deduce that \bqrn
   \ph_k(x)&=&\ph_j(x)^{\si}\,\ph_{S,\si}(\hat{T}_j(x))\\
   (\ph_j(x)^{\si})^{-1}\ph_k(x) &=& \ph_{S,\si}(\hat{T}_k(x))
\eqrn
Since the distribution of $(\ph_j(x)^{\si})^{-1}\ph_k(x)$
is purely atomic, it follows from Theorem~\ref{T:standard-quotients}\,(b)
that $[S]\in\OA^*(X,\GA)$ and
\[
    [S]=[\hat{T}_k\circ \hat{T}_j^{-1}]=[T_k][T_j]^{-1}\in\OA^*(X,\GA)
\]
contrary to the choice of $[T_i]$-s. Hence the induction step is
verified and the proof of Theorem~\ref{T:mainN} is completed. \qed

\subsection*{Proof of Corollary~\ref{C:entropy}}

Suppose that $[\Rout\Rel_{X,\GA}:\OA^*(X,\GA)]\ge n>1$.
Theorem~\ref{T:mainN} provides a $\GA$-equivariant quotient map
\[
   \pi:(X,\mu,\GA)\overto{} (G^{n-1}/\GA^{n-1},m_{G^{n-1}/\GA^{n-1}},\GA)
\]
where in the right hand side $\GA$ acts diagonally in each of the factors
$(G/\GA,m_{G/\GA},\GA^{\ta_i})$. For diagonal actions
the entropy is additive, so for every $\ga\in\GA$ one has
\bqrn
   h(X,\ga)&\ge& h(G^{n-1}/\GA^{n-1},\,m_{G^{n-1}/\GA^{n-1}},\,\ga)\\
   &=&\sum_{i=1}^{n-1} h(G/\GA,m_{G/\GA},\ga^{\ta_i})=(n-1)\cdot\ch(\Ad\ga)
\eqrn
which gives (\ref{e:entropy}).

In the context of smooth actions of $\GA$ on a compact
$d$-manifold $X$ another application of superrigidity for cocycles
allows to express the entropies $h(X,\mu,\ga)$ of elements
$\ga\in\GA$ via eigenvalues of $d$-dimensional
$G$-representations. More precisely, (see Furstenberg \cite{Furs}
Thm 8.3, or Zimmer \cite{Zi-book} 9.4.15) either $h(X,\mu,\ga)=0$
for all $\ga\in\GA$, or $h(X,\mu,\ga)=\ch(\rh(\ga))$, $\ga\in\GA$,
for some representation $\rh:G\to\GL_d(\bC)$. In particular one
has
\begin{equation}\label{e:ent-est}
     \inf_\ga \frac{h(X,\mu,\ga)}{\ch(\Ad\ga)}\le
   \max_{\dim\rh\le d}\inf_\ga\frac{\ch(\rh(\ga))}{\ch(\Ad\ga)}.
\end{equation}
Let us point out that in the above cited references the
$\GA$-action on $X$ and the measure $\mu$ were assumed to be
$C^2$-smooth, in order to apply Pesin's formula. However for the
\emph{inequality} (\ref{e:ent-est}) one only needs the upper bound
\[
  h(X,\mu,\ga)\le \max_{\dim\rh\le d} \ch(\rh(\ga)),\qquad\ga\in\GA
\]
which, being based on Margulis-Ruelle inequality, holds under
$C^1$-assumption on the action and does not require any regularity
assumptions on the measure $\mu$.

Using Borel's density theorem one may extend the $\inf$ in
(\ref{e:ent-est}) from $\ga\in\GA$ to $g\in G$, obtaining the
claimed estimate
\[
 	[\Rout\Rel_{X,\GA}:\OA^*(X,\GA)]\le 1+W_G(d).
\]
For a given $G$ the function $W_G(d)$ can be computed explicitly
in terms of the weights of irreducible representations, but here
let us confine the discussion to a general estimate $W_G(d)\le
d^2/8$, suggested to me by Dave Witte, whom I would like to thank.
For $k\ge 2$ let $\si_k$ denote the (unique !) irreducible
representation $\si_k$ of $\SL_2(\bR)$ in dimension $k$. If $h$
denotes the element ${\rm diag}(e,e^{-1})\in\SL_2(\bR)$, then the
eigenvalues of $\si_k(h)$ are $\{e^{k+1-2i} \mid i=1,\dots,k\}$ so
that
\[
   \ch(\si_k(h))=\sum_{i\le k/2} (k+1-2i)\le k^2/4
\]
Given a $d$-dimensional $G$-representation $\rh$ choose a subgroup
$\SL_2(\bR)\simeq G_0\subset G$, and let $g\in G$ correspond to
$h\in G_0$ above. The restriction $\rh|_{G_0}$ of $\rh$ to $G_0$
splits as a direct sum of irreducible $G_0$-representations
$\si_{d_i}$ with $\sum d_i=d$. Thus
\[
   \ch(\rh(g))=\sum \ch(\si_{d_i}(h))\le 1/4 \sum d_i^2\le d^2/4
\]
At the same time $\ch(\Ad_{G}(g))\ge\ch(\Ad_{\SL_2(\bR)}(h))=2$,
which gives
\[
   W_G(d)\le d^2/8.
\]
\qed

\section{Standard Examples without $G/\GA$ quotients} \label{S:examples}
In this section we prove \ref{T:torus}--\ref{C:ex-H-LA} applying
Theorem~\ref{T:main}.

\subsection*{Proof of Theorem~\ref{T:torus}}
Let us first verify the ergodicity and aperiodicity of the
$\GA$-action on $\bT^N$. Let $f\in L^2(\bT^N)\mapsto
\hat{f}\in\ell^2(\bZ^N)$ denote the Fourier transform. For
$A\in\SL_N(\bZ)$ one has $\widehat{f\circ A}=A^t\hat{f}$.
Therefore if $f\in L^2(\bT^N)$ is an invariant vector for a
subgroup $\LA\subset\SL_N(\bZ)$ then $\hat{f}\in\ell^2(\bZ^N)$ is
a $\LA^t$-invariant vector, and $\hat{f}$ is supported on finite
$\LA^t$-orbits on $\bZ^N$. Thus if $\GA$ fails to act ergodically
on $\bT^N$, then $\rh(\GA)^t$ has a non-trivial finite orbit on
$\bZ^N$, and for some finite index subgroup $\GA^\prime\subs\GA$
there is a non-trivial fixed vector for $\rh(\GA^\prime)^t$ in
$\bZ^N\subset\bR^N$. Since $\rh:G\to\SL_N(\bR)$ is rational,
Borel's density theorem implies that all of
$\rh(G)^t\subset\SL_N(\bR)$ has a non-trivial fixed vector, and
since $\rh(G)$ is totally reducible $\rh(G)$ also has non-trivial
fixed vectors contrary to the assumption. Thus $\GA$ acts
ergodically on $\bT^N$, and since the arguments apply to any
finite index subgroup of $\GA$, this action is aperiodic.

The $\GA$-action on $\bT^N$ can be assumed to be free. Indeed
$\SL_N(\bZ)$ acts freely (mod 0) on $\bT^N$ and so does
$\rh(\GA)\iso\GA$.

Next we claim that the system $(\bT^N,\GA)$ does not have $(\Ad
G/\GA^\prime,\Ad \GA)$ as a measurable quotient. In the case of
$\GA\subset\SL_n(\bZ)$ acting on $\bT^n$, $n>2$, this is easily
seen from the entropy comparison: for $\ga\in\SL_n(\bZ)$ with
eigenvalues $\la_1,\dots,\la_n$ one has
\[
  h(\bT^n,\ga)=\sum_i \log^+|\la_i|,\qquad
  h(\Ad G/\GA^\prime,\ga)=\sum_{i,j} \log^+|\la_i/\la_j|
\]
where $\GA^\prime$ is any lattice in $\Ad G=\PSL_n(\bR)$. Since
$|\det\ga|=1$, i.e. $\sum \log|\la_i|=0$, one has a strict
inequality $h(\bT^n,\ga)<h(\Ad G/\GA^\prime,\ga)$ as soon as $\ga$
has at least one eigenvalue off the unit circle. For the general
case we resort to a more complicated argument described below.

Now Theorem~\ref{T:main} (or rather its simple modification needed
to handle finite center) gives
\[
   \Rout\Rel_{\bT^N,\GA}\iso\OA^*(\bT^N,\GA)\iso\Aaut^*(\bT^N,\GA)/\GA.
\]
Evidently any $\si\in\GL_N(\bZ)$ which normalizes $\rh(\GA)$ gives
rise to the map $T_\si:x\mapsto \si(x)$ of $\bT^N$ which lies in
$\Aaut^*(\bT^N,\GA)$.
\begin{clm}\label{CLM:torus}
The correspondence $\si\to T_\si$ is an isomorphism
\[
   N_{\GL_N(\bZ)}(\rh(\GA))\iso\Aaut^*(\bT^N,\GA).
\]
\end{clm}
The correspondence $\si\to T_\si$ is clearly a monomorphism of
groups. To show its surjectivity consider a general
$T\in\Aaut^\ta(\bT^N,\GA)$ and let $\nu$ denote the lift of the
Lebesgue probability measure $m_{\bT^N}$ on $\bT^N$ to the graph
of $T$. Thus $\nu$ is a probability measure on
$\bT^N\times\bT^N=(\bR^N\times\bR^N)/(\bZ^N\times\bZ^N)$ which is
invariant and ergodic for the $(\rh\times\rh\circ\ta)(\GA)$-action
$\ga:(x,y)\mapsto (\rh(\ga)(x),\rh(\ga^\ta)(y))$. Witte's
\cite{Wi} Corollary~5.8 (based on Ratner's theorem) allows to
conclude that $\nu$ is a \emph{homogeneous} measure for some
closed subgroup
\[
   M\subs(\rh\times\rh\circ\ta)(\GA)\ltimes(\bR^N\times\bR^N)
\]
The connected component $M_0$ of the identity of $M$ can be viewed
as a subgroup of $\bR^N\times\bR^N$. The fact that $\nu$ is a lift
of $m_{\bT^N}$ to a graph of a m.p. bijection $T:\bT^N\to\bT^N$,
and the fact that $\bR^N$ is connected while $\bZ^N$ is discrete,
leads to the conclusion that $M_0\subset \bR^N\times \bR^N$
projects onto $\bR^N$ in both factors in a one-to-one fashion.
Hence $M_0=\{(x,\si(x)) \mid x\in \bR^N\}$ where $\si\in\Gaut
\bR^N$ which preserves $\bZ^N$, i.e. $\si\in\GL_N(\bZ)$, and $T$
has the form: $T(x)=\si(x)+t$ where $t\in\bT^N$ is such that
\[
    \si\circ\rh(\ga)(x)+t=\rh(\ga^\ta)(\si(x)+t).
\]
The latter means that $t$ is $\rh(\GA)$-fixed and
$\si\rh(\ga)\si^{-1}=\rh(\ga^\ta)$. An argument similar to the one
for aperiodicity of the action (based on the assumption that
$\rh(G)$ has no non-trivial fixed vectors), implies that $t$ has
to be trivial, so that $T$ is of the form $T_\si$ where $\si\in
N_{\GL_N(\bZ)}(\rh(\GA))$. The Claim is proved.

It remains to show that $\bT^N$ does not have $\Ad G/\GA^\prime$
as a measurable $\GA$-equivariant quotient. It follows from
Witte's \cite{Wi} 5.8 that measurable $\GA$-equivariant quotients
of $\bT^N=\bR^N/\bZ^N$ have the form $K\bs\bR^N/\LA$ where
$\bZ^N\subs\LA\subs\bR^N$ is a closed $\GA$-invariant subgroup and
$K$ is a closed subgroup of $\Aff(\bR^N/\LA)$ centralizing $\GA$;
moreover $K$ is acting non-ergodically on $\bR^N/\LA$. The latter space
can be identified with a quotient torus $\bT^n$, $n\le N$, on
which $\GA$ still acts by automorphisms, so that $K$ becomes a
subgroup of $\GL_n(\bZ)\ltimes\bT^n$. We claim that the
$\GA$-action on $K\bs\bT^n$ cannot be measurably isomorphic to the
$\GA$-action on $\Ad G/\GA^\prime$ because the former cannot be
extended to a $G$-action. In fact the $\GA$-action on $K\bs \bT^n$
cannot be extended to a measurable action of the smaller group -
the commensurator
\[
    \DE:={\rm Comm}_G(\GA)=\left\{g\in G\mid [\GA:g^{-1}\GA
    g\cap\GA]<\infty\right\}
\]
which is a dense subgroup in $G$ (follows from Margulis'
arithmeticity results \cite{Ma}). Indeed, let $g\mapsto T_g$,
$g\in\DE$, denote a hypothetical extension of the $\GA$-action on
$K\bs\bT^n$ to some measure-preserving $\DE$-action. For any
$g\in\DE$ there are finite index subgroups $\GA_1,\GA_2\subs\GA$
so that $\ta_g:\ga\mapsto g\ga g^{-1}$ is an isomorphism
$\GA_1\to\GA_2$. Thus $T_g$ satisfies $T_g(\ga\cdot x)=\ta_g(\ga)
\cdot T_g(x)$ for a.e. $x\in K\bs\bT^n$ and all $\ga\in\GA_1$.
Arguing as in the proof of Claim~\ref{CLM:torus} one shows that
such $T_g$ has to have an "algebraic" form, i.e. to be induced by
a linear map $\rh(g)\in\SL_n(\bR)$ which has to preserve the
lattice $\bZ^n$. The fact that the embedding $\GA\to\SL_n(\bZ)$
cannot be extended to the commensurator $\DE\supset\GA$ gives the
required contradiction. \qed

\subsection*{Proof of Theorem~\ref{T:ex-compact}}
By Margulis' Normal Subgroup Theorem \cite{Ma} (4.10) the
homomorphism $\rh:\GA\to K$ is actually an \emph{embedding}
(recall that $G$ and hence $\GA$ are assumed to be center free).
Thus without loss of generality we can assume that the proper
subgroup $L\subset K$ does not contain non-trivial normal factors
of $K$ (otherwise dividing by these factors we still remain in the
same setup). This means that the $K$-action $k_1:kL\to k_1kL$ is
free (mod 0) and so is the ergodic $\GA$-action
$(K/L,m_{K/L},\GA)$. This $\GA$-action is aperiodic: being
connected $K$ admits no proper closed subgroups of finite index,
and therefore any subgroup $\GA_1\subset \GA$ of finite index has
a dense image $\rh(\GA_1)$ in $K$ and acts ergodically on
$(K,m_K)$ as well as on its quotient $(K/L,m_{K/L})$. Furthermore,
such an action is irreducible - see Zimmer \cite{Zi-product} Prop
2.4. Clearly the discrete spectrum $\GA$-action on $K/L$ cannot
have equivarient quotients of the form $G/\GA$. Hence
Theorem~\ref{T:main} gives
\[
   \Rout\Rel_{(K/L,\GA)}=\OA^*(K/L,\GA)
\]

In Theorem~\ref{T:ex-compact} $K/L$ is a homogeneous space (recall
that being connected $K$ has to be a Lie group). However,
Theorem~\ref{T:ex-homogeneous} (or Ratner's theorem, in general)
does not apply to this situation because the acting group is not
generated by $\Ad$-unipotent elements. Yet the following general
result describing $\Aaut^*(K/L,\GA)$ can easily be obtained by
direct methods.
\begin{prop}\label{P:K-aut}
Let $K$ be a compact group, $\GA\subset K$ a dense subgroup and
$L\subs K$ a closed subgroup. Then the left $\GA$-action on
$(K/L,m_{K/L})$ is ergodic and
\bqrn
   \Aaut(K/L,\GA)&\iso& N_K(L)/L\\
   \Aaut^*(K/L,\GA)&\iso& N_{\Aff(K/L)}(\GA)
\eqrn
\end{prop}
\begin{rem}\label{R:Aut-K}
In the particular case of $L=\{e\}$ the first assertion, i.e. the
isomorphism $\Aaut(K,\GA)\iso K$, is easy seen as follows. Any
$T\in\Maut(K,m_K)$ can be written as $T(k)=k t_k^{-1}$ where
$k\mapsto t_k\in K$ is a measurable map. Then $T(\ga\cdot
k)=\ga\cdot T(k)$ translates into an a.e. identity $t_{\ga\cdot
k}=t_k$. Since $\GA$ acts ergodically on $(K,m_K)$ the map
$k\mapsto t_k$ is a.e. a constant $t\in K$, i.e. $T(k)=k t$. The
correspondence $T\in\Aaut(K,\GA)\mapsto t\in K$ is easily seen to
be an isomorphism of groups.
\end{rem}
\subsection*{Proof of Proposition~\ref{P:K-aut}}

Given $T\in\Aaut^\ta(K/L,\GA)$ let $\nu$ be the lift of $m_{K/L}$
to the graph of $T$ on $K/L\times K/L$, and let
\[
   R:=\{(k_1,k_2)\in K\times K \mid (k_1,k_2)_*\nu=\nu\}.
\]
$R\subs K\times K$ forms a closed (hence compact) group, containing
$\{(\ga,\ga^\ta) \mid \ga\in\GA\}$. The projections $p_i(R)$
of $R$ to $K$ are closed and contain $\GA$. Hence
$R$ projects onto $K$ in both coordinates. We claim that
\[
   R_1:=\{ k\in K \mid (k,e)\in R\},\qquad
   R_2:=\{ k\in K \mid (e,k)\in R\}
\]
are closed normal subgroups in $K$. Indeed, for $r_1\in R_1$ and
$k\in K$ there exists a $k_2\in K$ so that $(k,k_2)\in R$, and
\[
   (k,k_2)^{-1}(r_1,e)\,(k,k_2)=(k^{-1} r_1\,k,\,e)\in R
\]
shows that $k^{-1}r_1 k\in R_1$. Thus $R_1\normal K$ and similarly
$R_2\normal K$.

Since $\nu$ disintegrates into Dirac measures with respect to
$m_{K/L}$ under the projections $p_i:(K/L)\times (K/L)\to K/L$,
the $R_i$-actions on $K/L$ should fix $m_{K/L}$-a.e. point of
$K/L$. This means that $R_i\subs L$, and since $L$ is assumed not
to contain non-trivial normal factors of $K$, $R_i=\{e\}$ for
$i=1,2$. Hence $R$ has the form
\[
    R=\{ (k,\theta(k)) \mid k\in K \}
\]
for some bijection $\theta:K\to K$ which has to be a continuous
isomorphism, because $R\subset K\times K$ is a closed subgroup.

By definition of $R$ for all $k\in K$ and $m_{K/L}$-a.e. $k_1L$,
the point $(k k_1L, \theta(k)T(k_1L))$ is on the graph of $T$, i.e.
$T(kk_1L)=\theta(k)T(k_1L)$. Thus $T$ has the form $T(kL)=\theta(k)tL$
where $t\in K$ is such that $\theta(L)=tLt^{-1}$. Such $T$ can also
be written as $T(kL)=t\si(k)L$ where $\si(k)=t^{-1}\theta(k)t$, in
which case $\si\in N_{\Gaut K}(L)$. Thus $T$ comes from an
\emph{affine} map $a_{\si,t}\in\Aff(K/L)$. We conclude that
$\Aaut^*(K/L,\GA)$ coincides with $N_{\Aff(K/L)}(\GA)$.

Finally, an affine map $a_{\si,t}$ is in $\Aaut(K/L,\GA)$ if for
all $\ga\in\GA$ and a.e. $kL$
\[
   \ga t\si(k)L=t\si(\ga k)L=t\si(\ga)\si(k)L
\]
In view of the standing assumption that $L$ does not contain
normal subgroups of $K$ this means that $\si(\ga)=t^{-1}\ga t$ for
$\ga\in\GA$. Since $\GA$ is dense in $K$ we have $\si(k)=t^{-1} k
t$ for all $k\in K$ and $\si(L)=L$ means $t\in N_K(L)$. Hence
$a_{\si,t}:kL\mapsto (tt^{-1})ktL=ktL$ and $t,t^\prime\in N_K(L)$
give rise to the same map of $\Aff(K/L)$ iff $t^\prime t^{-1}\in
L$. This gives the desired identification
\[
    \Aaut(K/L,\GA)\iso N_K(L)/L
\]
\qed
This completes the proof of Theorem~\ref{T:ex-compact}.

\subsection*{Proof of Theorem~\ref{T:ex-homogeneous}}
By Theorem~\ref{T:quotients} the system $(H/\LA,m_{H/\LA},\GA)$
has a $\GA$-equivariant quotient map
\[
   \pi:(H/\LA,m_{H/\LA})\overto{} (G/\GA^\prime,m_{G/\GA^\prime})
\]
only if there exists a surjective continuous homomorphism
$\si:H\to G$ with $\si(\LA)\subs \GA^\prime\iso\GA$ and
$\si\circ\rh(\ga)=t\ga t^{-1}$ for some $t\in G$. An existence of
such a homomorphism $\si$ was explicitly excluded by the
assumption, so that Theorem~\ref{T:main} gives
$\Rout\Rel_{H/\LA,\GA}=\OA^(H/\LA,\GA)\iso\Aaut^*(H/\LA,\GA)/\GA$.
To identify $\Aaut^*(H/\LA,\GA)$ we invoke
Theorem~\ref{T:quotients} again to conclude that
\[
    \Aaut^*(H/\LA,\GA)= N_{\Aff(H/\LA)}(\rh(\GA))
\]
which presents $\OA^*(H/\LA,\GA)$ as the quotient of
$N_{\Aff(H/\LA)}(\rh(\GA))$ by the image of
$\GA\overto{\rh}H\hookrightarrow\Aff(H/\LA).$ One also has
$\Aaut(H/\LA,\GA)\iso C_{\Aff(H/\LA)}(\rh(\GA))$. \qed

\subsection*{Proof of Corollary~\ref{C:ex-H-LA}}
If $\rh:G\to H$ is an embedding (or isomorphism) of $G$ into
another semi-simple real Lie group $H$ (center free and without
compact factors) and $\LA\subset H$ is an irreducible lattice,
then the $G$-action on $(H/\LA,m_{H/\LA})$ is free and by
Howe-Moore's theorem is not only ergodic but actually mixing.
Hence also the restriction of this action to $\GA$-action is free
and mixing, and in particular irreducible and aperiodic. The
assumptions of the Corollary guarantee that there does not exits
an epimorphism $\si:H\to G$ with $\si(\LA)\subs\GA$, so that
Theorem~\ref{T:ex-homogeneous} applies showing
\[
   \Rout\Rel_{H/\LA,\GA}=\OA^*(H/\LA,\GA)\iso
   \Aaut^*(H/\LA,\GA)/\GA\iso N_{\Aff(H/\LA)}(\rh(\GA))/\rh(\GA).
\]
Recal that $\Aff(H/\LA)$ contains $H$ as a subgroup of finite
index dividing $|\Gout\LA|$. 
Hence, upon passing to a subgroup of index
dividing $|\Gout\LA|$, the group $\Rout\Rel_{H/\LA,\GA}\iso
N_{\Aff(H/\LA)}(\rh(\GA))/\rh(\GA)$ can be reduced to
$N_H(\rh(\GA))/\rh(\GA)$, which contains the centralizer
$C_H(\rh(\GA))=C_H(\rh(G))$ as a subgroup of index dividing
$|\Gout\GA|$. \qed

\section{Proof of Theorem~\ref{T:G-mod-GA}} \label{S:with-G-GA-factor}
\subsection*{Case $(G/\GA,\GA)$}
Choose a two-sided fundamental domain $X\subset G$ for $\GA$ and
define the transformation $I:X\to X$ by $I:x\mapsto x^{-1}\GA\cap
X$. Note that both $X$ and $X^{-1}$ are two-sided, in particular
right, fundamental domains and therefore $I$ is a measurable
bijection of $X$. Moreover,
\[
   I(\ga\cdot x)=I(\ga x \la(\ga,x)^{-1})=\la(\ga,x) x^{-1}\GA\cap X
   =\la(\ga,x)\cdot I(x)
\]
which means that $I\in\Raut\Rel_{(G/\GA,\GA)}$ and the corresponding
rearrangement cocycle $\al_I$ is $\la=\la_X:\GA\times X\to\GA$. Observe that
\[
   \ga\cdot x=\ga x \la(\ga,x)^{-1}\qquad{\rm means\ that}\qquad
   \la(\ga,x)=(\ga\cdot x)^{-1}\,\ga\, x
\]
(with the usual multiplication in $G$ on the right hand side), so that
the embedding $X\to G$ is precisely the ``straightening map'' $\ph$ corresponding
to the cocycle $\al_I=\la_X$ and the trivial automorphism $\ta_0:\ga\mapsto \ga$;
in other words $\ph_{I,\ta_0}(x)=x$. 
From Theorem~\ref{T:standard-quotients}\,(a1)
we conclude that $[I]\not\in\OA^*(G/\GA,\GA)$ and therefore
\begin{equation}\label{e:2}
   [\Rout\Rel_{(G/\GA,\GA)}:\OA^*(G/\GA,\GA)]\ge 2
\end{equation}
while Theorem~\ref{T:mainN} (or Corollary \ref{C:entropy}) show
that this index is at most two proving an equality in (\ref{e:2}).
Theorem~\ref{T:quotients} gives
\[
   \Aaut^*(G/\GA,\GA)\iso N_{\Aff(G/\GA)}(\GA)
\]
Note that an affine map $a_{\si,t}\in\Aff(G/\GA)$
($a_{\si,t}:g\GA\mapsto t\si(g)\GA$ where $\si\in N_{\Gaut
G}(\GA)$ and $t\in G$) satisfies
\[
    a_{\si,t}(\ga\cdot g\GA)=\ga^\ta\cdot a_{\si,t}(g\GA)
\]
iff $\si(\ga)=t^{-1}\ta(\ga)t$, in particular $t\in N_G(\GA)$.
Thus $\Aaut^*(G/\GA,\GA)\iso N_{\Aff(G/\GA)}(\GA)\iso N_{\Gaut
G}(\GA)$, with $g\GA\mapsto g^\ta\GA$, $\ta\in N_{\Gaut
G}(\GA)\iso\Gaut\GA$, giving all twisted action automorphisms.
Hence $\OA^*(G/\GA,\GA)\iso \Gaut \GA/\GA\iso\Gout\GA$. Since this
group commutes with $[I]$, we obtain
\[
   \Rout\Rel_{(G/\GA,\GA)}\iso\bZ/2\bZ\times\Gout(\GA)
\]
as claimed.

Before turning to the systems $(G^n/\GA^n,\GA)$ for general finite $n\ge 1$,
observe that $G/\GA$ can be viewed as the factor of
$G^2_e:=\{(g,g^{-1})\in G\times G \mid g\in G\}$ modulo the relation
$(g,g^{-1})\sim (g\ga_1,\ga_1^{-1}g)$, $\ga_1\in\GA$. With this identification
$G/\GA\iso (G^2_e/\sim)$ the left $\GA$-action on $G/\GA$ corresponds to
the quotient of the action
$\ga: (g,g^{-1})\mapsto (\ga g,g^{-1}\ga^{-1})$ modulo $\sim$,
while the map $I$ arises from the flip $(g,g^{-1})\mapsto(g^{-1},g)$.

\subsection*{Case $(G^n/\GA^n,\GA)$, $n\in\bN$.}
Given a general finite $n$ consider the  set
\[
   G^{n+1}_e:=\{ (g_0,\dots,g_n)\in G_e^{n+1} \mid g_0\cdots g_n=e\}
\]
with the natural measure and an equivalence $\sim$ defined by
\[
   (g_0,g_1,\dots,g_{n-1},g_n)\ \sim\
   (g_0\ga_1^{-1},\,\ga_1 g_1\ga_2^{-1},\dots,\ga_{n-1}g_n\ga_n^{-1},\, \ga_n g_n)
\]
for $\ga_1,\dots,\ga_n\in\GA$. The map $p:G_e^{n+1}\to
(G/\GA)^n=G^n/\GA^n$ given by
\[
   p:(g_0,\dots,g_n)\mapsto (g_0\GA,\,g_0g_1\GA,\,\dots,\,g_0 g_1\cdots g_{n-1}\GA)
\]
factors through a bijection $q:(G_e^{n+1}/\sim)\,\to G^n/\GA^n$.
Note that the following $\GA$-action on $G_e^{n+1}$
\[
   \ga\,:\,(g_0,\,g_1,\,\dots,\,g_{n-1},\,g_n)\mapsto
   (\ga g_0,\,g_1,\,\dots,\,g_{n-1},\,g_n\ga^{-1})
\]
descends to an action on $(G_e^{n+1}/\sim)$ which is isomorphic, via $q$, to the
diagonal $\GA$-action on $G^n/\GA^n$
\[
   \ga:(g_1\GA,\dots,g_n\GA)\,\mapsto\, (\ga g_1\GA,\,\dots,\,\ga g_n\GA).
\]
The cyclic permutation $\tilde{T}$ of order $(n+1)$
\[
   \tilde{T}\,:\,(g_0,g_1,\dots,g_{n-1},g_{n})\mapsto (g_1,g_2,\dots,g_n,g_0)
\]
is easily seen to preserve the $\GA$-orbits on
$(G_e^{n+1}/\sim)\iso G^n/\GA^n$ and thereby defines a relation
automorphism $T\in\Raut\Rel_{G^n/\GA^n,\GA}$ with
$[T^{n+1}]\in\OA(G^n/\GA^n,\GA)$.

We would like to present $T$ as an explicit transformation of
$(G^n/\GA^n,m_{G^n/\GA^n})$ as follows. The cocycle
$\la_X:\GA\times X\to\GA$ corresponding to the two-sided
fundamental domain $X\subset G$ can be extended to a cocycle of
$G$, i.e. $\la=\la_X:G\times X\to\GA$ still defined by $gx \in
\,X\,\la(g,x)$. The left $G$-action on $X\iso G/\GA$ can thus be
written as
\[
    g\cdot x=g\,x\,\la(g,x)^{-1}
\]
where on the right hand side we use the usual multiplication in $G$.
Using these notations and viewing $x\in X\subset G$ both as points of
the space $X$ and as $G$-elements one obtains an explicit form for $T$:
\[
   T\,:\,(x_1,\dots,x_n)\mapsto
         (x_1^{-1}\cdot x_2,\,x_1^{-1}\cdot x_3,\,\dots,
          \,x_1^{-1}\cdot x_n,\,I(x_1)).
\]
Observe that
\bqrn
   &&T(\ga\cdot(x_1,\dots,x_n))=
       T(\ga x_1\la(\ga,x_1)^{-1},\,\dots,\,\ga x_n\la(\ga,x_n)^{-1})\\
    &&= (\la(\ga,x_1)x_1^{-1}\cdot x_2,\,\la(\ga,x_1) x_1^{-1}\cdot x_3,\,\dots,\,
         \la(\ga,x_1)\cdot I(x_1))\\
    &&= \la(\ga,x_1)\cdot T(x_1,\dots,x_n).
\eqrn Hence $T\in\Raut\Rel_{(G^n/\GA^n,\GA)}$ with the
rearrangement cocycle being
\[
   \al_{T}(\ga,(x_1,\dots,x_n))=\la(\ga,x_1).
\]
A similar computation shows that for $1\le k\le n$ one has
\[
   \al_{T^k}(\ga,(x_1,\dots,x_n))=\la(\ga,x_k)
\]
and therefore the corresponding ``straightening'' map is given by
\[
   \ph_{T^k,\ta_0}(x_1,\dots,x_n)=x_k\in G.
\]
It now follows from Theorem \ref{T:standard-quotients}\,(a1) that
$T^k\not\in\OA^*(G^n/\GA^n,\GA)$ for $k=1,\dots,n$. In particular
\[
   [\Rout\Rel_{(G^n/\GA^n,\GA)}:\OA^*(G^n/\GA^n,\GA)]\ge n+1
\]
which is, in fact, an equality due to the upper bound $(n+1)$
provided by Theorem~\ref{T:mainN} (or Corollary~\ref{C:entropy}).

To identify $\OA^*(G^n/\GA^n,\GA)$ we invoke the second part of
Theorem~\ref{T:quotients} with $H:=G^n$ and $\LA:=\GA^n$ and note
that affine maps of $G^n/\GA^n$ have the form
\[
   (g_1\GA,\dots,g_n\GA)\,\mapsto\,
   (t_1g_{p(1)}^{\ta_1}\GA,\dots,t_ng_{p(n)}^{\ta_n}\GA)
\]
where $p\in S_n$ is a permutation of $\{1,\dots,n\}$, $\ta_i\in
N_{\Gaut G}(\GA)\iso\Gaut\GA$ and $t_i\in G$. One easily checks
that such a map normalizes the diagonal $\GA$-action iff
$\ta_1=\dots=\ta_n=\ta$ and $t_1=\dots=t_n=t$ where $t\in
N_G(\GA)$. Hence $\Aaut^*(G^n/\GA^n,\,\GA)$ consists of the maps
\[
   S_{p,\ta}:(g_1\GA,\dots,g_n\GA)\mapsto
   (g_{p(1)}^\ta\GA,\dots,g_{p(n)}^\ta\GA)
\]
where $p\in S_n$ and $\ta\in N_{\Gaut G}(\GA))$. The obvious
relation $S_{p,\ta}\circ
S_{p^\prime,\ta^\prime}=S_{pp^\prime,\ta\ta^\prime}$ gives
$\Aaut^*(G^n/\GA^n,\GA)\iso S_n\times N_{\Gaut}(\GA)$ and
\[
   \OA^*(G^n/\GA^n,\,\GA)\iso
    S_n\times (N_{\Gaut}(\GA)/\GA)\iso S_n\times\Gout(\GA).
\]
$\Rout\Rel_{(G^n/\GA^n,\GA)}$ is generated by $[T]$ and $\OA^*(G^n/\GA^n,\GA)$,
and the explicit form of $T$ and $S_{p,\ta}$ allows one to check that
\[
   \Rout\Rel_{(G^n/\GA^n,\GA)}\iso S_{n+1}\times\Gout(\GA)
\]
as claimed.

\subsection*{Case  $(G^\infty/\GA^\infty,\GA)$}
Finally, let us turn to the case of $n=\infty$, i.e. the diagonal $\GA$-action
on $(\ol{X},\ol{\mu}):=(G/\GA,m_{G/\GA})^{\bZ}$.
Choose a two-sided fundamental domain
$X\subset G$, so that $\ol{X}=X^{\bZ}$, and let
$\la=\la_X:G\times X\to \GA$ and $I:X\to X$ be as before.
Consider the map $T:\ol{X}\to \ol{X}$ defined by
\[
   T:(\dots,x_{-1},x_0,x_1,\dots)\mapsto
     (\dots,\,x_1^{-1}\cdot x_0,\,I(x_1),\,x_1^{-1}\cdot x_2,\dots)
\]
so that for $k\neq 0$
\[
   (T^k \bar{x})_i:=\left\{\begin{array}{lll}
     x_{k}^{-1}\cdot x_{i+k}&\qquad &i\neq 1-k\\
     I(x_{k}) &\qquad& i=1-k
   \end{array}\right.
\]
and observe that
\[
   T^k(\ga\cdot\bar{x})=\la(\ga,x_k)\cdot T^k(\bar{x}).
\]
As before, for $k\neq 0$ we have $\al_{T^k}(\ga,\bar{x})=\la(\ga,x_k)$
and $\ph_{T^k,\ta_o}(\bar{x})=x_k$ so that $[T]^k\not\in\OA^*(\ol{X},\GA)$.
\begin{clm}\label{CLM:T-generates}
$\Rout\Rel_{\ol{X},\GA}$ is generated by $[T]$ and $\OA^*(\ol{X},\GA)$.
\end{clm}
(Note that in previous cases similar statement followed immediately
from the upper bound provided by Corollary~\ref{C:entropy}).
Choose an $S\in\Raut\Rel_{\ol{X},\GA}\setminus\OA^*(\ol{X},\GA)$ and let
\[
   \pi:\ol{X}\to G/\GA,\qquad\pi_*\ol{\mu}=m_{G/\GA},\qquad
   \pi(\ga\cdot\bar{x})=\ga^\ta\cdot\pi(\bar{x})
\]
be the standard quotient map provided by Theorem~\ref{T:standard-quotients}.
\begin{lem}\label{L:pi-xk}
$\pi(\bar{x})=x_k^\ta$ for some $k\in\bZ$ and $\ta\in N_{\Gaut G}(\GA)$.
\end{lem}
\Pf
Denote by $\ol{\nu}$ the probability measure on $(G/\GA)^{\bZ}\times (G/\GA)$
obtained by the lift of $\ol{\mu}$ to the graph of $\pi$. Fix an $r\in\bN$,
let
\[
   H:=\prod_{-r}^r G\qquad \LA:=\prod_{-r}^r \GA
\]
and let $p:G^{\bZ}\to H$ be the projection on $\{-r,\dots,r\}$-coordinates.
Denote by $\ol{\nu}^{(r)}$ the $p\times \Id$-projection of $\ol{\nu}$ to
$H/\LA\times G/\GA$. Then one can deduce from Theorem~\ref{T:joinings} that
either
\begin{itemize}
\item[(i)] $\ol{\nu}^{(r)}=m_{H/\LA}\times m_{G/\GA}$, or
\item[(ii)] There exists $k\in\{-r,\dots,r\}$, $\ta\in N_{\Gaut G}(\GA)$
  so that for any $F\in C_c(H/\LA\times G/\GA)$
  \[
     \int F\,d\ol{\nu}^{(r)}
     =\int F(x_1,\dots,x_k,\dots,x_n,x_k^\ta)\,
     dm_{G/\GA}(x_1)\cdots dm_{G/\GA}(x_n).
  \]
\end{itemize}
As $r\to\infty$ case (i) cannot persist forever,
because that would imply that $\ol{\nu}=\ol{\mu}\times m_{G/\GA}$
which is impossible. On the other hand as soon as (ii) occurs,
the index $k$ and $\ta\in N_{\Gaut G}(\GA)$ do not change.
This proves the Lemma.
\qed

With the explicit form of $\pi:\ol{X}\to G/\GA$ provided by
Lemma~\ref{L:pi-xk} we invoke Theorem~\ref{T:standard-quotients}\,(a3)
to conclude that there exists $\hat{S}\in\Raut\Rel_{\ol{X},\GA}$
with $[S]=[\hat{S}]$, $\ta\in N_{\Gaut G}(\GA)$ and $k\neq0\in\bZ$
so that
\[
    \ph_{\hat{S},\ta}(\ol{x})\GA=(x_k)^\ta\GA.
\]
Recalling that also for $T^k$ we have
$\ph_{T^k,\ta_0}(\ol{x})=x_k$ one concludes that
$[S]=[\hat{S}]\in[T^k]\OA^*(\ol{X},\GA)$ using the same argument
as in the proof of Theorem~\ref{T:mainN}. This completes the proof
of Claim~\ref{CLM:T-generates}. \qed

Any permutation $p$ of $\bZ$ and any $\ta\in N_{\Gaut G}(\GA)$
give rise to the map $S_{p,\ta}\in\Aaut^\ta(\ol{X},\GA)$
\[
   S_{p,\ta}\,:\,
   (g_i\GA)_{i\in\bZ}\,\mapsto\,(g_{p(i)}^\ta\GA)_{i\in\bZ}.
\]
On the other hand if $S\in\Aaut^*(\ol{X},\GA)$ let $\ol{\nu}$ on
$\ol{X}\times\ol{X}$ be the lift of $\ol{\mu}$ to the graph of $S$
and let $\ol{\nu}^r$ be the projection of this measure to
$\prod_{-r}^r G/\GA\times\prod_{-r}^r G/\GA$. Then applying the
Joining Theorem~\ref{T:joinings} to this \emph{finite dimensional}
situation successively for $r\to\infty$, one concludes that such
$S$ has to be of the form $S_{p,\ta}$. Hence \bqrn
   \Aaut^*\left(\ol{X},\,\GA\right)&\iso& S_\infty\times\ N_{\Gaut G}(\GA)\\
   \OA^*\left(\ol{X},\,\GA\right)&\iso& S_\infty\times\ \Gout\GA.
\eqrn
Finally, the explicit form of $[T]$ and $[S_{p,\ta}]$ allows
to conclude that
\[
   \Rout\Rel{(\ol{X},\GA)}\iso S_{\infty+1}\times\Gout\GA
\]
where the symbols $S_\infty$ and $S_{\infty+1}$ can be interpreted as
the inclusion of the permutation
group of $\bZ$ in the permutation group of $\bZ\cup\{pt\}$.
\qed
\section{Proof of Theorem~\ref{T:SLnZp}} \label{S:SLnZp}

Throught this section $\GA=\PSL_n(\bZ)$, $G=\PSL_n(\bR)$ and $n\ge3$.
Let $S_0=\{p_1,\dots,p_r\}$ be a given finite set of primes and consider the ergodic 
$\GA$-action on the compact profinite group $K=\prod_{p\in S_0}\PSL_n(\bZ_p)$.
We denote $H=\prod_{p\in S_0}\PSL_n(\bQ_p)$ and $\LA=\PSL_n(\bZ[S^{-1}])\subset H$. 
Then $\LA$ is a dense countable subgroup of locally compact totally disconnected group
$H$ and $\GA=\LA\cap K$.

Following Gefter \cite{Ge2} we first observe that $\Rout\Rel_{K,\GA}$ contains $H$.
Indeed restricting the type ${\rm II}_\infty$ relation $\Rel_{H,\LA}$ to $K$ we obtain a
type ${\rm II}_1$ relation $\Rel_{K,\GA}=\Rel_{H,\LA}\cap({K\times K})$ and  
\[
	\Rout\Rel_{K,\GA}\iso\Rout\Rel_{H,\LA}\supseteq \OA(H,\LA)\iso H
\]
using the straight forward ${\rm II}_\infty$-type generalizations of Lemmas \ref{L:induction}, 
\ref{L:ICC}(a) and the remark following \ref{P:K-aut} respectively. 

We need to find explicit representatives $T_h\in\Raut\Rel_{K,\GA}$ for $h\in H$,
so that $h\mapsto[T_h]$ is the above imbedding.
Since $K$ is open and $\LA$ is dense in $H$, given any $h\in H$, there exist $\la_0\in\LA$ and $k_0\in K$ 
so that $h=\la_0k_0$.
The maps
\[
	\tilde{T}_h: x\mapsto xh,\rmand\tilde{T}^\prime_h:x\mapsto \la_0^{-1} xh,\qquad(x\in H)
\]
are in $\Raut\Rel_{H,\LA}$ and $[\tilde{T}_h]=[\tilde{T}^\prime_h]\in\Rout\Rel_{H,\LA}$.
Denoting  the open compact subgroup $\la_0K\la_0^{-1}\cap K$ by $K_1$, note that
\[
	\tilde{T}^\prime_h(K_1)\subset K\quad{\rm because}\quad 
	\tilde{T}^\prime_h(x)=\la_0^{-1} x\la_0 k_0\in (\la_0^{-1}K_1\la_0) k_0\subset K
\]
Thus for $x,y\in K_1$ we have 
\[
	(x,y)\in\Rel_{K,\GA}=\Rel_{H,\LA}|_{K}\qquad{\rm iff}\qquad
	(\tilde{T}^\prime_h(x),\tilde{T}^\prime_h(y))\in\Rel_{H,\LA}|_{K}=\Rel_{K,\GA}.
\] 
Therefore $\tilde{T}^\prime_h|_{K_1}$ is a restriction of some $T_h\in\Raut\Rel_{K,\GA}$,
with its outer class $[T_h]\in\Rout\Rel_{K,\GA}$ being uniquely defined by the initial $h\in H$. 
Denoting $\GA_1:=\la_0\GA\la_0^{-1}\cap\GA$ a finite index subgroup of $\GA$ which is
dense in $K_1$, we observe that the restriction of the rearrangement cocycle $\al_{T_h}$ to 
$\GA_1\times K_1$ is
\begin{equation}\label{e:Th}
	\al_{T_h}(\ga_1,x_1)=\la_0^{-1}\ga_1\la_0\qquad(\ga_1\in\GA_1,\ {\rm a.e.\ }x_1\in K_1)
\end{equation}
The automorphism $\ga_1\mapsto \la_0^{-1}\ga_1\la_0$ of the lattice $\GA_1\subset G$ 
extends to an inner (given by $\la_0\in G$) automorphism of $G$, so in terms of the Standard Quotients Theorem~\ref{T:standard-quotients} the class $[\ta]\in\Gout G$ associated to such 
$[T_h]\in\Rout\Rel_{K,\GA}$ is always trivial.
On the other hand the transpose map $T_0:(k_1,\dots,k_r)\mapsto (k_1^t,\dots,k_r^t)$
which is clearly in $\Raut\Rel_{K,\GA}$ defines the unique outer element $[\ta]\in\Rout G$ 
(take $\ta(g)=(g^t)^{-1}$). 
One easily checks that the group generated by $[T_0]$ and $[T_h]$, $h\in H$, in $\Rout\Rel_{K,\GA}$
is $\bZ/2$-extension of $H$.

We shall now prove that the latter group is all of $\Rout\Rel_{K,\GA}$. 
Take any $[T]\in\Rout\Rel_{K,\GA}$. 
Possibly composing with $T_0$ we may assume that $[\ta]\in\Rout G\iso \bZ/2$ associated
with $[T]$ is trivial, and will show that such $[T]$ is $[T_h]$ for some $h\in H$.
Applying the Standard Quotients Theorem we may take $\ta$ to be the identity on $G$.
Since $(K,\GA)$ cannot have $(G/\GA,m_{G/\GA},\GA)$ among its measurable quotients, 
we deduce  that
\begin{enumerate}
\item 
   There exists a finite $\GA$-orbit $F=\{g_1\GA,\dots,g_k\GA\}\subset G/\GA$,
   and a measurable $\GA$-equivariant map $\pi:K\to F$ with 
   \[
   	\pi(\ga x)=\ga \pi(x)\qquad (\ga\in\GA,\ x\in K)
   \]
\item   
   Let $\GA_i=\GA\cap g_i\GA g_i^{-1}$, $i=1,\dots,k$ - these are conjugate
   subgroups of index $k$ in $\GA$; 
   the sets $X_i=\pi^{-1}(\{g_i\GA\})\subset K$ are $\GA_i$-invariant and ergodic
   measurable subsets with $\mu(X_i)=1/k$; if $K_i$ is the closure of $\GA_i$ in $K$
   then $X_i=K_i y_i$ (mod 0) -- cosets of $K_i$-s; as the latter are open and compact
   subsets of $K$ we obtain an open partition into disjoint sets which we still denote by 
   $X_i$. Up to reordering we may assume that $X_1$ contains the identity of $K$,
   i.e. $X_1=K_1$.
\item
   There exists $\hat{T}\in\Raut\Rel_{K,\GA}$ with $[\hat{T}]=[T]\in\Rout\Rel_{K,\GA}$ so that 
   \[
   	\al_{\hat{T}}(\ga_1,x_1)=g_1^{-1} \ga_1 g_1\qquad (\ga_1\in\GA_1,\ x_1\in X_1).
   \]
\end{enumerate}
Note that the last formula resembles (\ref{e:Th}).
Property (1) means that $g_1\in {\rm Comm}_G(\GA)=\PSL_n(\bQ)$.
\begin{clm}\label{CLM:g1inLA}
$g_1\in\LA=\PSL_n(\bZ[S^{-1}])$.
\end{clm}
\Pf
Let us expand the notations slightly:  for an arbitrary finite set $S$ of primes let
\[
	K_S=\prod_{p\in S}\PSL_n(\bZ_p),\qquad \LA_S=\PSL_n(\bZ[S^{-1}])
\]
and let $\mu_S$ denote the normalized Haar measure on $K_S$.
We shall denote by $\ol{\GA_1}^S$ the closure of the index $k$ subgroup
$\GA_1=\GA\cap g_1\GA g_1^{-1}\subset \GA$ in $K_S$.
The $\GA_1$-ergodic component $X_1\subset K=K_{S_0}$ is a coset
of the open compact subgroup $\ol{\GA_1}^{S_0}$ of $K$ and by (2)
\[
	\frac{1}{k}=\mu(X_1)=\mu_{S_0}(\ol{\GA_1}^{S_0}).
\]
Let $S_1$ be the set of primes appearing in the denominators of $g_1\in\PSL_n(\bQ)$,
i.e. $S_1$ is the smallest set of primes (possibly empty) such that $g_1\in\LA_{S_1}$.

It follows from the Strong Approximation Theorem that if $S=S^\prime\sqcup S^\pprime$ 
is a disjoint union of two finite sets of primes,  then 
\[
	\ol{\GA_1}^S=\ol{\GA_1}^{S^\prime}\times\ol{\GA_1}^{S^\pprime}\subs 
	K_{S^\prime}\times K_{S^\pprime}=K_S,
\]
and  
$\ol{\GA_1}^{S^\pprime}=K_{S^\pprime}$ if and only if $S^\prime\cap S_1=\emptyset$.
On the other hand if $S_1\subset S^\prime$ then it is easy to see that
\[
	\mu_{S^\prime}(\ol{\GA_1}^{S^\prime})=\frac{1}{[\GA:\GA_1]}=\frac{1}{k}.
\]
Writing $S=S_0\cup S_1=S_0\sqcup S_2$ where $S_2=S_1\setminus S_0$ we have
\[
	\frac{1}{k}=\mu_{S_0}(\ol{\GA_1}^{S_0})\ge 
		\mu_{S_0}(\ol{\GA_1}^{S_0})\cdot \mu_{S_2}(\ol{\GA_1}^{S_2})
		= \mu_{S}(\ol{\GA_1}^{S})=\frac{1}{[\GA:\GA_1]}=\frac{1}{k}.
\]
So $\mu_{S_2}(\ol{\GA_1}^{S_2})=1$, that is $\ol{\GA_1}^{S_2}=K_{S_2}$,
which means that $S_2=\emptyset$ and $S_1\subs S_0$ as claimed.
\qed

Having proved that $g_1\in\LA$, we recall that by (3) the original $T\in\Raut\Rel_{K,\GA}$
can be replaced by $\hat{T}$ with $[T]=[\hat{T}]\in\Rout\Rel_{K,\GA}$ so that
\begin{equation}\label{e:That}
	\hat{T}(\ga_1 x_1)=g_1^{-1}\ga_1 g_1\,\hat{T}(x_1)
\end{equation}
for all $\ga_1\in\GA_1$ and $\mu$-a.e. $x_1\in X_1$.
We have also made sure that $X_1=K_1$ - the closure
of $\GA_1=\GA\cap g_1\GA g_1^{-1}$ in $K$.
\begin{clm}
$\ \hat{T}(k)=g_1^{-1}k g_1 z_1$ for some fixed $z_1\in K$ and a.e. $k\in K_1$.
\end{clm}
\Pf
The map $\ga_1\mapsto g_1^{-1}\ga_1g_1$ is an isomorphism between 
finite index subgroups $\GA_1\to\GA_1^\prime:= g_1^{-1}\GA g_1\cap\GA$ of $\GA$.
It extends to an isomorphism $K_1\to K_1^\prime$ between open compact subgroups of $K$,
where $K_1^\prime$ is the closure of $\GA_1^\prime$ in $K$. 
(Note that $K_1=K\cap g_1 K g_1^{-1}$ and $K_1^\prime=g_1^{-1} K g_1\cap K$ as subsets of $H$). 

Let $X_1^\prime=\hat{T}(X_1)\subset K$. In view of (\ref{e:That}), $X_1^\prime$ is
one of the $\GA_1^\prime$-ergodic components of $X_1$, and therefore is a single
$K_1^\prime$-coset, $X_1^\prime=K_1^\prime y$ for some $y\in X_1^\prime$.
Let $R:K_1\to K_1$ be the composition of the following maps
\[
	K_1=X_1\to X_1^\prime\to K_1^\prime\to K_1\qquad R(k)=g_1 \hat{T}(k)y^{-1} g_1^{-1}.
\]
In view of (\ref{e:That}) we have for all $\ga\in\GA_1$ and $\mu$-a.e. $k\in K_1$:
\[
	R(\ga k)=g_1 g_1^{-1}\ga g_1\hat{T}(k)y^{-1} g_1^{-1}=\ga R(k). 
\]  
Since $\GA_1$ is dense in the compact group $K_1$, we have $R(k)=k k_0$ 
for some fixed $k_0\in K_1$ and a.e. $k\in K_1$ (see Proposition~\ref{P:K-aut} and the following Remark).
This allows us to compute 
\[
	\hat{T}(k)=g_1^{-1}k k_0 g_1 y=g_1^{-1}k g_1 z_1\rmwhere z_1=(g_1^{-1} k_0 g_1)y\in K.
\]
\qed

Taking $h=g_1 z_1\in H$ we observe that the map $T_h\in\Raut\Rel_{K,\GA}$, discussed in the first
part of this section, agrees with $\hat{T}$ on a positive measure subset $K_1\subset K$,
and therefore (as in the proof of Lemma~\ref{L:induction}) 
\[
	[T]=[\hat{T}]=[T_h]\in\Rout\Rel_{K,\GA}
\]
which completes the proof of the Theorem. 
\qed


\end{document}